\documentclass{amsart}
\usepackage{latexsym}
\usepackage{amsmath, amssymb, enumerate}
\usepackage[all]{xy}
\newtheorem{dfs}{Definition}[section]
\newtheorem{lms}[dfs]{Lemma}
\newtheorem{thms}[dfs]{Theorem}
\newtheorem{props}[dfs]{Proposition}
\newtheorem{qus}[dfs]{Question}
\newtheorem{cors}[dfs]{Corollary}
\newtheorem{rems}[dfs]{Remark}

\newtheorem{nota}[dfs]{Notation}
\newtheorem{blank}[dfs]{}
\newcommand{\Z}{\ensuremath{\mathcal{Z}}}
\newcommand{\dt}{\ensuremath{d_{\tau}}}

\newcommand{\wt}{\ensuremath{\widetilde{W}(A)}}
\newcommand{\co}{\ensuremath{\mathrm{C}}}
\newcommand{\ma}{\ensuremath{\mathrm{M}}}
\newcommand{\aff}{\ensuremath{\mathrm{Aff}}}
\newcommand{\laff}{\ensuremath{\mathrm{LAff}}}

\author{Francesc Perera and Andrew S. Toms}
\address{Departament de Matem\`atiques, Universitat Aut\`onoma de Barcelona, 08193 Bellaterra, Barcelona, Spain}
\email{perera@mat.uab.es}
\address{Department of Mathematics, York University, 4700 Keele St., Toronto, ON, Canada, M3J 1P3}
\email{atoms@mathstat.yorku.ca}
\thanks{Research supported by the DGI MEC-FEDER through Project MTM2005-00934, and the
Comissionat per Universitats i Recerca de la Generalitat de
Catalunya. The second named author was also supported in part by an NSERC Discovery Grant}

\title{Recasting the Elliott conjecture}
\date{\today}

\begin{document}

\begin{abstract}
Let $A$ be a simple, unital, finite, and exact C$^*$-algebra
which absorbs the Jiang-Su algebra $\mathcal{Z}$ tensorially.  We prove that the Cuntz semigroup
of $A$ admits a complete order
embedding into an ordered semigroup which is obtained from the Elliott invariant in a
functorial manner.  We conjecture that this embedding is an isomorphism, and prove the
conjecture in several cases.  In these same cases --- $\mathcal{Z}$-stable algebras all
--- we prove that the Elliott conjecture in its strongest form is equivalent to a
conjecture which appears much weaker.  Outside the class of $\mathcal{Z}$-stable C$^*$-algebras,
this weaker conjecture has no known counterexamples, and
it is plausible that none exist.  Thus, we reconcile the still intact principle of Elliott's classification
conjecture --- that $\mathrm{K}$-theoretic invariants will classify separable and nuclear C$^*$-algebras
--- with the recent appearance of counterexamples to its strongest concrete form.
\end{abstract}

\maketitle
\section{Introduction}
The Elliott conjecture for C$^*$-algebras operates on
two levels:  on the one hand, it is a meta-conjecture asserting that separable
and nuclear C$^*$-algebras will be classified up to $*$-isomorphism by
$\mathrm{K}$-theoretic invariants;  on the other, it is a collection of
concrete classification conjectures, where the $\mathrm{K}$-theoretic invariants
in question are specified and depend on the class of algebras being considered.
In the case of stable Kirchberg algebras (simple, nuclear, purely infinite, and satisfying
the Universal Coefficients Theorem), the correct invariant is the graded Abelian
group $\mathrm{K}_0 \oplus \mathrm{K}_1$ (\cite{K}, \cite{P}).  For non-simple algebras of real rank zero, $\mathrm{K}$-theory with
coefficients seems to suffice (\cite{DG}, \cite{Ei}).  For a unital, separable, and
nuclear C$^*$-algebra $A$, the invariant
\[
I(A) := \left(
(\mathrm{K}_0(A),\mathrm{K}_0(A)^+,[1_A]),\mathrm{K}_1(A),\mathrm{T}(A),
r_A \right)
\]
--- topological $\mathrm{K}$-theory, the (possibly empty) Choquet simplex $\mathrm{T}(A)$ of tracial states,
and the pairing $r_A: \mathrm{T}(A) \times \mathrm{K}_0(A) \to
\mathbb{R}$ given by evaluating a trace at a $\mathrm{K}_0$-class
--- is known as the Elliott invariant, and has been very
successful in confirming Elliott's conjecture for simple algebras.

In its most general form, the Elliott conjecture may be stated as follows:
\begin{blank}
[Elliott, c. 1989]\label{ec}
There is a $\mathrm{K}$-theoretic functor $F$ from the category of separable and nuclear
C$^*$-algebras such that if $A$ and $B$ are separable and nuclear, and there is an
isomorphism
\[
\phi: F(A) \to F(B),
\]
then there is a $*$-isomorphism
\[
\Phi: A \to B
\]
such that $F(\Phi) = \phi$.
\end{blank}

\noindent
We will let (EC) denote the conjecture above with the Elliott invariant $I(\bullet)$
substituted for $F(\bullet)$, and with the class of algebras under consideration
restricted to those which are simple and unital.  (EC) has been shown
to hold in many situations. An exhaustive list of these results would be impossibly
long, but \cite{El1}, \cite{El3}, \cite{EG}, \cite{EGL}, \cite{K}, and \cite{Li1} are among
the most important works.
We refer the reader to R{\o}rdam's book (\cite{R3}) for a comprehensive
overview of Elliott's classification programme.

Recent examples due first to R{\o}rdam and later the second named author have shown the
currently proposed invariants (i.e., the proposed values of $F$ in Conjecture \ref{ec}) to be
insufficient for the classification of all simple, separable, and nuclear
C$^*$-algebras (\cite{R2}, \cite{T1}, \cite{T2}).  In particular, (EC) does
not hold.  There are two options:  enlarge the
proposed invariants, or restrict the class of algebras considered.

The Cuntz semigroup of a C$^*$-algebra $A$ is a positively ordered
Abelian semigroup whose elements are equivalence classes of
positive elements in matrix algebras over $A$ (see section 2 for
details).  Let $W(A)$ denote this semigroup, and let $\langle a
\rangle$ denote the equivalence class of a positive element $a \in
\mathrm{M}_n(A)$. The semigroup $W(A)$ may be thought of as a
generalisation of the semigroup $V(A)$ of Murray-von Neumann
equivalence classes of projections in matrices over $A$, provided
that $A$ is stably finite. Theorem 1 of \cite{T2} states that
there exist simple, separable, nuclear, and non-isomorphic
C$^*$-algebras which agree on each continuous and homotopy
invariant functor from the category of C$^*$-algebras, and which
furthermore have the same simplex of tracial states.  These algebras are
distinguished by their Cuntz semigroups, whence this invariant is
extremely sensitive. (Indeed, it is already unmanageably large for
commutative C$^*$-algebras with contractible spectrum --- see
\cite[Lemma 5.1]{T2}.) It thus suggests itself as the minimum
quantity by which the Elliott invariant $I(\bullet)$ ought to be
enlarged.  The sequel will be concerned in large part with the
relationship between (EC) and the following statement:
\begin{blank}[WEC]
Let $A$ and $B$ be simple, separable, unital, and nuclear
C$^*$-algebras.  If there is an isomorphism
\[
\phi\colon \left(W(A), \langle 1_A \rangle, I(A)\right) \to
\left(W(B), \langle 1_B \rangle, I(B)\right),
\]
then there is a $*$-isomorphism $\Phi\colon A \to B$ which induces $\phi$.
\end{blank}
\noindent
There are no known counterexamples to the conjecture (WEC) among stably finite
algebras, and perhaps none
exist.  But asking for the Cuntz semigroup as part of the invariant seems
strong indeed, given its sensitivity and the fact that (EC) alone is so often true.
The theme of the sequel is that (WEC)
and (EC) are reconciled upon restriction to the largest class of C$^*$-algebras for which
(EC) may be expected to hold.  (WEC) may thus be viewed as the appropriate
specification of the Elliott conjecture for simple, separable, unital, nuclear, and stably
finite C$^*$-algebras.
(We have, for the time being,
glossed over what exactly is meant by isomorphism at the level of invariants in
both (EC) and (WEC), so as not to burden this introduction with technicalities.
The appropriate notions of isomorphism will be introduced in section 4.)

It is generally agreed that the largest restricted class of
algebras for which (EC) can hold consists of those algebras which
absorb the Jiang-Su algebra $\mathcal{Z}$ tensorially
(\cite{JS1}).  Indeed, this fact is obvious if one considers only
algebras with weakly unperforated ordered $\mathrm{K}_0$-groups (a
condition which holds in every confirmation of (EC)) --- by
Theorem 1 of \cite{gjsu}, the tensor product of such an algebra,
say $A$, with $\mathcal{Z}$ has the same Elliott invariant as $A$,
and so (EC) predicts that $A \cong A \otimes \mathcal{Z}$.  If $A$
is any C$^*$-algebra and the minimal tensor product $A \otimes
\mathcal{Z}$ is isomorphic to $A$, then we say that $A$ is {\it
$\mathcal{Z}$-stable}. Our first main result is:

\begin{thms}\label{ecwec1}
Upon restriction to $\mathcal{Z}$-stable C$^*$-algebras, (EC) implies (WEC).
\end{thms}

\noindent
Notice that this theorem does not follow from the mere fact that the invariant
considered in (WEC) is finer that the Elliott invariant.  This is due to the functorial
nature of Elliott-type conjectures:  an isomorphism at the level of the invariant must
lift to an isomorphism at the level of C$^*$-algebras which, moreover, induces the original
isomorphism of invariants.

More surprising, perhaps, is this:

\begin{thms}\label{equiv}
Let $\mathcal{C}$ denote the class of all simple, unital, separable, nuclear,
and $\mathcal{Z}$-stable C$^*$-algebras $A$ which are either
\begin{enumerate}
\item[{\rm (i)}] of real rank zero, or
\item[{\rm (ii)}] have finitely many pure tracial states.
\end{enumerate}
Then, (EC) and (WEC) are equivalent in $\mathcal{C}$.
Moreover, there is a functor $G$ from the category of Elliott invariants to the category of Elliott
invariants augmented by the Cuntz semigroup such that
\[
G(I(A)) = \left(W(A), \langle 1_A \rangle, I(A)\right).
\]
\end{thms}
\noindent
In proving Theorem \ref{equiv} we shall see that an algebra $A \in \mathcal{C}$ has, 
up to Cuntz equivalence, relatively few positive elements.  This contrasts sharply with
the counterexample to (EC) in \cite{T2}.  Significant is the
fact that $A$ need not be of real rank zero;  it may be projectionless but for zero and
the unit.  Most progress on (EC) from a general point of view has so far required the real
rank zero assumption.  We also outline a proof that Theorem \ref{equiv} holds for Goodearl
algebras, so that conditions (i) and (ii) of the theorem are, in principle, removeable.  (Indeed,
we conjecture as much.)  The
proof of Theorem \ref{equiv} gives the first calculations of Cuntz semigroups for C$^*$-algebras without
the real rank zero property, and even in the real rank zero case generalises considerably
the earlier results of Blackadar and Handelman (\cite{bh}).

The paper is organised as follows:  in Section 2 we recall the definition of the Cuntz semigroup,
and establish several results about its order structure;  in Section 3 we compute $W(\mathcal{Z})$,
and examine the basic structure of $W(A \otimes \mathcal{Z})$; Section 4 contains an embedding
theorem which establishes Theorem \ref{ecwec1}; Section 5 contains a calculation of the
Grothendieck enveloping group of the Cuntz semigroup for finite $\mathcal{Z}$-stable algebras;
Sections 6 and 7 are devoted to proving Theorem \ref{equiv} in cases (i) and (ii),
respectively; in Section 8 we sketch a proof of Theorem \ref{equiv} for Goodearl algebras;
Section 9 raises some questions for future research.

\vspace{2mm}
\noindent
\emph{Acknowledgements.} The second named author would like to thank George Elliott
for several inspiring conversations related to the results herein, for his hospitality
at the Fields Institute in the summer of 2005, and for his guidance and support in general.
Some of the work on this paper was carried out during the second named author's visit
to the Centre de Recerca Matem\`atica at the Universitat Aut\`onoma de Barcelona, hosted
by the first named author and Pere Ara.  The support of all concerned with this visit
is gratefully acknowledged.  Finally, we thank the referee, whose careful reading and suggestions
improved our exposition considerably.

\section{The Cuntz semigroup and comparison}

J. Cuntz introduced in~\cite{Cu} a notion of comparison between positive elements in a C$^*$-algebra
that extends the usual (Murray-von Neumann) comparison for projections.  This allowed him to prove the
existence of dimension functions in stably finite simple C$^*$-algebras.
(The assumption of simplicity was subsequently removed by D. Handelman in~\cite{han}.)

Explicitly, if $a$ and $b$ are positive elements in a C$^*$-algebra $A$, then we write $a\precsim b$
provided there is a sequence of elements $(x_n)$ in $A$ such that $a=\lim\limits_{n\to \infty}x_n bx_n^*$.
This relation can be extended to the (local) C$^*$-algebra $M_{\infty}(A)$ defined as the inductive limit
of $M_n(A)$ via the inclusion mappings $M_n(A)\hookrightarrow M_{n+1}(A)$ given by
$x\mapsto\left(\begin{smallmatrix} x & 0 \\
 0 & 0\end{smallmatrix}\right)$.
Let $\mathrm{M}_{\infty}(A)_+$ denote the set of positive elements in $\mathrm{M}_{\infty}(A)$.
For elements $a$, $b$
in $M_{\infty}(A)_+$, we write $a\precsim b$ provided that $a\precsim b$ in $M_n(A)$ for some $n$ such
that $a$, $b\in M_n(A)$. (If we view $a$ and $b$ in two different sized matrices over $A$, the above is
equivalent to having $a=\lim\limits_{n\to\infty} x_nbx_n^*$ where the $x_n$ are suitable rectangular matrices.)
If both $a\precsim b$ and $b\precsim a$, we will write $a\sim b$ and call $a$ and $b$ \emph{Cuntz equivalent}.
We shall denote the equivalence class of an element $a$ in $M_{\infty}(A)_+$ by $\langle a\rangle$, and we
will in this paper denote the set of all such equivalence classes by $W(A)$ (although this notation is not
uniform in the literature). For $a, b \in M_{\infty}(A)_+$ we write $a\oplus b$ for
the element $\left(\begin{smallmatrix} a & 0 \\ 0 & b \end{smallmatrix}\right) \in M_{\infty}(A)_+$.
If $\langle a\rangle$, $\langle b\rangle\in W(A)$, we define $\langle a\rangle+\langle b\rangle=\langle a\oplus b\rangle$.
It is easy to verify that this operation is does not depend on the representatives chosen and that $W(A)$ becomes an
Abelian semigroup with identity element $\langle 0\rangle$ (and thus an Abelian monoid). We shall refer
to $W(A)$ as the \emph{Cuntz semigroup of $A$}. All semigroups in this paper will be Abelian and assumed to have an identity element,
which we shall denote by $0$.

Recall that projections $p,q \in M_{\infty}(A)$ are
\emph{Murray-von Neumann equivalent} ($p\sim q$) if there is an
element $x$ in $M_{\infty}(A)$ such that $p=xx^*$ and $q=x^*x$;
$p$ is \emph{subequivalent} to $q$ (in symbols $p\precsim q$) if
there is a projection $q' \in \mathrm{M}_{\infty}(A)$ such that
$p\sim q'$ and $q'\leq q$.  The notions of Murray-von
Neumann equivalence and Cuntz equivalence coincide
for the set of projections in matrices over a stably finite
C$^*$-algebra, but do not coincide in general.
Let $[p]$ denote the Murray-von Neumann equivalence
class of $p$. The set of all such equivalence classes is denoted
$V(A)$, and is also an Abelian semigroup (with identity element
$[0]$) under the operation $[p]+[q]=[p\oplus q]$. There is a
natural semigroup morphism $\varphi\colon V(A)\to W(A)$, given by
$[p]\mapsto \langle p\rangle$, which is injective if $A$ is stably
finite. In this case, we identify $V(A)$ with its image under
$\varphi$.

\begin{dfs}\label{purepos}
Let $A$ be a C$^*$-algebra, and let $W(A)_+$ denote the subset of $W(A)$ consisting of classes which are not
the classes of projections.  If $a\in A_+$ and $\langle a\rangle\in W(A)_+$, then we
will say that $a$ is purely positive and denote the set of such elements by $A_{++}$.
\end{dfs}

One of the advantages of the relation $\precsim$ is that it allows the decomposition of elements up to arbitrary approximations.
If $\epsilon>0$ and $a\in A_+$, then $(a-\epsilon)_+$ will denote the positive part of $a-\epsilon\cdot 1$, that
is, $(a-\epsilon)_+=f(a)$, where $f\colon\mathbb{R}\to \mathbb{R}$ is given by $f(t)=\max\{t-\epsilon,0\}$. It is
proved in~\cite[Proposition 2.4]{Rfunct} (see also~\cite[Proposition 2.6]{KR}) that $a\precsim b$ if and only if
for any $\epsilon>0$, there exists $\delta>0$ and $x$ in $A$ such that $(a-\epsilon)_+=x(b-\delta)_+x^*$. (This is
in turn equivalent to the statement that, for any $\epsilon>0$, there is $\delta>0$ such that $(a-\epsilon)_+\precsim (b-\delta)_+$.)

The next proposition shows that despite the typically non-algebraic ordering
on the Cuntz semigroup, one can always complement projections.

\begin{props}\label{projcomplement}
Let $A$ be a C$^*$-algebra.  Let $a,p \in \mathrm{M}_{\infty}(A)_+$
be such that $p$ is a projection and $p \precsim a$. Then, there exists $b \in \mathrm{M}_{\infty}(A)_+$
such that $p \oplus b \sim a$.
\end{props}

\begin{proof}
By passing to a suitable matrix over $A$, we may assume that actually $p$, $a\in A$. Let $0<\epsilon<1$. Since
$p\precsim a$, we have that $p \sim (p-\epsilon)_+=xax^*$, for some $x$ in $pA$. Set $p'=a^{\frac{1}{2}}x^*xa^{\frac{1}{2}}$.
Then $p'$ is a projection equivalent to $p$ and $p'\leq \Vert x\Vert^2 a$, which is Cuntz equivalent to $a$.
Therefore we may assume at the outset that $p\leq a$.

We claim now that $p\oplus (1-p)a(1-p)\sim a$. By~\cite[Lemma 2.8]{KR}, we always have that $a\precsim pap\oplus
(1-p)a(1-p)$. Since $pap\leq \Vert a\Vert^2p\sim p$, we obtain that $a\precsim p\oplus (1-p)a(1-p)$. To establish
the converse subequivalence, it will suffice to show that both $p$ and $(1-p)a(1-p)$ belong to the hereditary
algebra $A_a$ generated by $a$, because then $p+(1-p)a(1-p)\in A_a$. From this it follows that $p+(1-p)a(1-p)\precsim a$.

By our assumption we have that $p\leq a$ and thus $p\in A_a$.
Also, $(1-p)a^{\frac{1}{2}}=a^{\frac{1}{2}}-pa^{\frac{1}{2}}\in A_a$, whence $(1-p)a(1-p)\in A_a$.
\end{proof}

Let $M$ be a preordered Abelian semigroup, with order relation
denoted by $\leq$. Recall that a non-zero element $u$ in $M$ is
said to be an \emph{order-unit} provided that for any $x$ in $M$
there is a natural number $n$ such that $x\leq nu$. A \emph{state}
on a preordered monoid $M$ with order-unit $u$ is an order
preserving monoid morphism $s\colon M\to \mathbb{R}$ such that
$s(u)=1$. We denote the (convex) set of states by
$\mathrm{S}(M,u)$. In the case of a unital C$^*$-algebra $A$, the
set of states on the Cuntz monoid $W(A)$ is referred as to the
\emph{dimension functions} on $A$ and denoted by $\mathrm{DF}(A)$
(see also~\cite{bh},~\cite{Rfunct},~\cite{P1}).

A dimension function $s$ is lower semicontinuous if $s(\langle
a\rangle)\leq \liminf\limits_{n\to\infty} s(\langle a_n\rangle)$
whenever $a_n\to a$ in norm. The set of all lower semicontinous
dimension functions on $A$ is denoted by $\mathrm{LDF}(A)$. Note
that any dimension function $s$ induces a function $d_s\colon
M_{\infty}(A)\to \mathbb{R}$ given by $d_s(a)=s\langle
a^*a\rangle$. With this notation, lower semicontinuity of $s$ as
defined above is equivalent to lower semicontinuity of the
function $d_s$.

We shall denote by $\mathrm{T}(A)$ the simplex of normalised traces defined
on a unital C$^*$-algebra $A$, and by $\mathrm{QT}(A)$ the
simplex of quasitraces. (We will work mostly with simple
unital C$^*$-algebras in the sequel, and so take the term ``quasitrace''
to mean a normalised 2-quasitrace --- see \cite{bh}.)  We have
$\mathrm{T}(A)\subseteq \mathrm{QT}(A)$, and equality holds if $A$
is exact and unital by the main theorem of~\cite{Ha}. Any quasitrace
$\tau$ defines a lower semicontinuous dimension function
\[
\dt(a)= \lim_{n \to \infty} \tau(a^{1/n}),
\]
provided that the domain of $\dt$ is restricted to positive elements.
In fact, it was proved in~\cite[Theorem II.2.2]{bh} that if $d\in
\mathrm{LDF}(A)$, then there is a unique quasitrace $\tau$ such that $d=\dt$.
It is clear that if $a\precsim b$, then for any dimension function $d$ we have $d(a)\leq d(b)$.

The next definition is not new (see~\cite{blsur},~\cite{P1}), but bears repeating.
\begin{dfs}\label{fcq+}
Let $A$ be a unital C$^*$-algebra with $a,b \in A_+ \backslash \{0\}$.
We say that $A$ has the Fundamental Comparability Property for Positive Elements,
denoted by (FCQ+), if $a \precsim b$ whenever $\dt(a) < \dt(b)$ for every
$\tau \in \mathrm{QT}(A)$. 
\end{dfs}
Villadsen gave the first example
of a simple C$^*$-algebra for which (FCQ+) fails (\cite{V1}).  In his example the positive
elements $a$ and $b$ are projections. (FCQ+) may hold for all pairs of
projections, yet fail in general (\cite{T2}). The abbreviation (FCQ+) derives from 
Blackadar's Fundamental Comparability Question, which asks if (FCQ+) holds whenever
$a$ and $b$ are projections.  In the literature, however, a C$^*$-algebra $A$ with
(FCQ+) is usually said to have \emph{strict comparison of positive elements} or
simply \emph{strict comparison}.  The latter terminology will be employed in the sequel.

\begin{lms}
\label{comparisonlemma}
Let $A$ be a unital C$^*$-algebra with $a\in A_{+}$. For any faithful quasitrace $\tau$
and $\epsilon,\eta,\delta \in \sigma(a)$ with $\epsilon<\eta<\delta$ we have $\dt ((a-\delta)_+)<\dt ((a-\epsilon)_+)$.
\end{lms}

\begin{proof}
Since $(a-\epsilon)_+$ and $(a-\delta)_+$ belong to the C$^*$-algebra $C^*(a)$ generated by
$a$, we may assume that $A=C^*(a)$. Then $\tau$ corresponds to a probability measure $\mu_\tau$ 
on $\sigma (a)$ which is nonzero on every open set. 
By~\cite[Proposition I.2.1]{bh} we have $\dt (b)=\mu_\tau(\mathrm{Coz}(b))$, where
$\mathrm{Coz}(b)$ is the cozero set of a nonnegative function $b$ in $C^*(a)$.

Put $U_\epsilon = \{ (\epsilon,\infty] \cap \sigma(a)\}$;  define $U_{\delta}$ similarly.
Let $V$ be an open subset of $\sigma(a)$ containing $\eta$ and such that $V \subseteq U_{\epsilon} \cap U_{\delta}^c$.
Let $b$ be a nonnegative function on $\sigma(a)$ such that $\mathrm{Coz}(b) = V$ and $b \leq (a-\epsilon)_+$.
Now $b$ is orthogonal to $(a-\delta)_+$ and $(a-\delta)_+ + b \leq (a-\epsilon)_+$, so
\[
\dt((a-\delta)_+) + \dt(b) \leq \dt((a-\epsilon)_+);
\]
$\dt(b) = \mu_{\tau}(V)$ is nonzero, and the lemma follows.
\end{proof}

\begin{rems}\label{specminfty} {\rm
We will occasionally refer to the spectrum $\sigma(a)$ of a positive element $a \in \mathrm{M}_{\infty}(A)$.
Since $a$ may be viewed as an element of arbitrarily large matrix algebras over $A$, we always assume that
$0 \in \sigma(a)$ for consistency. }
\end{rems}

\begin{props}\label{pureposcomparison}
Let $A$ be a simple C$^*$-algebra with strict comparison of positive elements.
Let $a \in A_{++}$ and $b \in A_+$ satisfy $\dt(a) \leq \dt(b)$ for every $\tau \in \mathrm{QT}(A)$.
Then, $a \precsim b$.
\end{props}

\begin{proof}
If $A$ has no quasitrace, then strict comparison of positive elements reduces to
the condition that for any nonzero positive elements $a,b \in A$, there is a sequence
$x_j$ in $A$ such that $x_j b x_j^* \to a$ as $j \to \infty$.  Thus, $a \precsim b$.

Suppose that $\mathrm{QT}(A)$ is nonempty.  Each quasitrace is faithful since $A$ is simple.
Since $a \in A_{++}$, we have that $a\neq 0$ and $0 \in \sigma(a)$.
Then, there is a strictly decreasing sequence
$\epsilon_n$ of positive reals in $\sigma(a)$ converging to zero.
We also know by~\cite[Section 6]{blsur} (see
also~\cite[Proposition 2.6]{KR}) that the set $\{x\in A_+\mid
x\precsim b\}$ is closed, and since $(a-\epsilon_n)_+\to a$ in
norm it suffices to prove that $(a-\epsilon_n)_+ \precsim b$ for
every $n \in \mathbb{N}$.

Let $\tau\in \mathrm{QT}(A)$ be given, and apply Lemma \ref{comparisonlemma}
with $\epsilon = 0$, $\eta=\epsilon_{n+1}$, and $\delta = \epsilon_n$ to see that
\[
\dt((a-\epsilon_n)_+)< \dt(a) \leq \dt(b).
\]
Using strict comparison we conclude that $(a-\epsilon_n)_+\precsim b$ for
all $n$, as desired.
\end{proof}

\begin{props}\label{projpureposcomp}
Let $A$ be as in Proposition~\ref{pureposcomparison}.  Let $p$ be a projection in $A$, and
let $a \in A_{++}$.  Then, $p \precsim a$ if and only if $\dt(p) < \dt(a)$ for each
$\tau \in \mathrm{QT}(A)$.
\end{props}

\begin{proof}
If $A$ has no quasitrace, then it is purely infinite and $p \precsim a$ (\cite{KR}).

Assume that $\mathrm{QT}(A)$ is nonempty.  The reverse implication follows from strict
comparison.  We prove the contrapositive of the forward implication.
Suppose that $\dt(a)\leq \dt(p)$ for some $\tau \in \mathrm{QT}(A)$, and
let $1 > \epsilon > 0$ be given. By~\cite[Proposition 2.4]{Rfunct} there exists a
$\delta > 0$ such that
\[
(p-\epsilon)_+ \precsim (a-\delta)_+.
\]
This implies that
\[
\dt((p-\epsilon)_+) \leq \dt((a-\delta)_+).
\]
But $p$ is a projection, so the functional calculus implies that
\[
\dt((p-\epsilon)_+) = \dt(p).
\]
Now
\[
\dt((p-\epsilon)_+) \leq \dt((a-\delta)_+) < \dt(a) \leq \dt(p) = \dt((p-\epsilon)_+),
\]
a contradiction.
\end{proof}

The hypotheses of Propositions \ref{pureposcomparison} and
\ref{projpureposcomp} are satisfied whenever $A$ is simple, unital,
and $W(A)$ satisfies the technical condition of being almost unperforated (see \cite{R1}).
In particular, $A$ could be a simple, unital and finite C$^*$-algebra
absorbing the Jiang-Su algebra $\Z$ tensorially (\cite[Corollary 4.6]{R1}).

\begin{props}\label{pureposspec}
Let $A$ be a simple, unital, and stably finite C$^*$-algebra, and let $a \in
\mathrm{M}_{\infty}(A)_+$.  Then, $\langle a \rangle = \langle p \rangle$ for
a projection $p$ in $\mathrm{M}_{\infty}(A)_+$ if and only if $0 \notin \sigma(a)$
or $0$ is an isolated point of $\sigma(a)$.
\end{props}

\begin{proof}
If $0 \notin \sigma(a)$ or $0$ is an isolated point of $\sigma(a)$, then
$\langle a \rangle = \langle p \rangle$ for a projection $p$ in $\mathrm{M}_{\infty}(A)_+$
by a straightforward functional calculus argument.

Now suppose that $\langle a \rangle = \langle p \rangle$ for a projection $p$ in
$\mathrm{M}_{\infty}(A)_+$ and $0$ is an accumulation point of $\sigma(a)$.
Choose $\epsilon \in [0,1) \cap \sigma(a)$ and a (necessarily faithful) quasitrace
$\tau \in \mathrm{QT}(A)$. Using~\cite[Proposition 2.4]{Rfunct},
there is $0<\delta \in \sigma(a)$ such that $(p-\epsilon)_+\precsim (a-\delta)_+$. If $\delta>\epsilon$, then
\[
\dt(a)=\dt(p)=\dt((p-\epsilon_+)\leq\dt((a-\delta)_+)\leq\dt
((a-\epsilon)_+)<\dt(a)\,
\]
by Lemma~\ref{comparisonlemma}; this is impossible. Thus
$\delta\leq\epsilon$, and by assumption we may find
$\delta'\in\sigma(a)$ such that $\delta'<\delta$. A second
application of Lemma~\ref{comparisonlemma} implies that
\[
\dt(a)=\dt(p)\leq\dt((a-\delta)_+)\leq
\dt((a-\delta')_+)<\dt(a)\,;
\]
this, too, is impossible.
%
\end{proof}

\noindent
If one replaces the assumptions of simplicity and being stably finite with stable rank one,
then Proposition \ref{pureposspec} is due to P. Ara (\cite[Proposition 3.12]{P1}).

\begin{cors}\label{varecov}
Let $A$ be a unital C$^*$-algebra which is either simple and stably finite or of stable rank one.  Then:
\begin{enumerate}
\item[(i)] 
$W(A)_+$ is a semigroup, and is absorbing in the sense that if one has $a \in W(A)$
and $b \in W(A)_+$, then $a + b \in W(A)_+$; 
\item[(ii)] $V(A)$ can be identified with the set of all $x \in W(A)$ which satisfy the
following condition: if $x\leq y$ for some $y\in W(A)$, then $x+z=y$ for some $z\in W(A)$.
\end{enumerate}
\end{cors}

\begin{proof}
For (i), take $\langle a \rangle, \langle b \rangle \in W(A)_+$ and notice that
the spectrum of $a \oplus b$ contains the union of the spectra of $a$ and
$b$.  Apply Proposition \ref{pureposspec}. 

For (ii), let $X$ be the following set:
\[
\{x\in W(A)\mid \text{ if }x\leq y \text{ for }y\in W(A)\,,\,\text{then }x+z=y\text{ for some }z\in W(A)\}.
\]
By Proposition~\ref{projcomplement}, we already know that $V(A)\subseteq X$.

Conversely, if $\langle x\rangle\in X$, then we may find a projection $p$ (in $M_{\infty}(A)$) such that
$\langle x\rangle\leq\langle p\rangle$. But then there is $z$ in $M_{\infty}(A)$ for which $x\oplus z\sim p$.
Since $0$ is an isolated point in $\sigma(p)$, the same will be true of
$\sigma (x)$. Invoking Proposition \ref{pureposspec} or \cite[Proposition 3.12]{P1} as appropriate,
we find a projection $q$ such that $q\sim x$, and so $\langle x\rangle\in V(A)$.
\end{proof}

The last proposition of this section, though straightforward, will be quite
important in the sequel.

\begin{props}\label{lsc}
Let $A$ be a stably finite unital C$^*$-algebra, and let $a \in A_+$.  Then,
the map $\tau \mapsto \dt(a)$ is a lower semicontinuous bounded function on
$\mathrm{T}(A)$.
\end{props}

\begin{proof}
Since $\langle \lambda a \rangle = \langle a \rangle$ for every $\lambda \in \mathbb{R}^+ \backslash \{0\}$,
we may assume that $||a|| \leq 1$.  Then, $f_n(\tau) := \tau(a^{1/n})$ is an increasing sequence
of continuous functions on $\mathrm{T}(A)$ with pointwise limit $f(\tau) := \dt(a)$.
\end{proof}




\section{$\Z$-stable C$^*$-algebras}

In this section we give a precise description of $W(\Z)$
(Theorem~\ref{WZ} below), and establish the important fact that
$W(\bullet)_+$ is a $\mathbb{R}^+$-cone for certain finite and
$\Z$-stable C$^*$-algebras.  In the study of the Cuntz semigroup
for simple, unital, and $\Z$-stable
C$^*$-algebras, the finite case is the only interesting one.
Indeed, a simple, unital, and $\Z$-stable C$^*$-algebra $A$ either
has stable rank one or is purely infinite (see~\cite[Theorem
3]{gjsu} and also~\cite[Corollary 5.1 and Theorem 6.7]{R1}). If
$A$ is purely infinite, then $a\precsim b$ for all non-zero
positive elements (see~\cite{LZ}).  It follows that
$W(A)=\{0,\langle 1 \rangle\}$ ($\langle 1 \rangle + \langle 1
\rangle = \langle 1 \rangle$), and that the Grothendieck group
$\mathrm{K}_0^*(A)$ of $W(A)$ is zero.

We begin with some notation.  For a compact convex set $K$, denote
by $\aff(K)^+$ the semigroup of all positive, affine, continuous,
and real-valued functions on $K$; $\laff (K)^+ \subseteq
\aff(K)^+$ is the subsemigroup of lower semicontinuous functions,
and $\laff_b(K)^+ \subseteq \laff (K)^+$ is the subsemigroup
consisting of those functions which are bounded above.  The use of
an additional ``+'' superscript (e.g., $\aff(K)^{++}$) indicates
that we are considering only strictly positive functions together
with the zero function.   Unless otherwise noted, the order on
these semigroups will be pointwise. $\aff (K)^+$ is algebraically
ordered with this ordering, but $\laff (K)^+$, in general, is not
(unless $K$ is, for example, finite dimensional).

Given two partially ordered semigroups $M$ and $N$, a
homomorphism $\varphi\colon M\to N$ is said to be an
\emph{order-embedding} provided that $\varphi(x)\leq \varphi(y)$
if and only if $x\leq y$. A surjective order-embedding will be
called an \emph{order-isomorphism}.

Let $\leq_{\mathbb{R}}$ denote the
usual order on the real numbers. We equip the disjoint union
$\mathbb{Z}^+\sqcup\mathbb{R}^{++}$ with a semigroup structure
by using the usual addition inside the components $\mathbb{Z}^+$
and $\mathbb{R}^{++}$ and declaring that $x+y\in\mathbb{R}^{++}$
whenever $x\in\mathbb{Z}^+$ and $y\in\mathbb{R}^{++}$.  Define
an order $\leq_{\Z}$ on this semigroup by using the usual order
inside the components $\mathbb{Z}^+$ and $\mathbb{R}^{++}$, and
the following order for comparing $x \in \mathbb{Z}^+$ and
$y \in \mathbb{R}^{++}$:  $x \leq_{\Z} y$ iff $x
<_{\mathbb{R}} y$, while $x \geq_{\Z} y$ iff $x \geq_{\mathbb{R}} y$.
With this ordering, $1_{\mathbb{Z}^+}$ is an order-unit.

\begin{thms}\label{WZ}
The ordered semigroup $(W(\Z),\langle 1_{\Z} \rangle)$ is
order-isomorphic (as an ordered monoid with order-unit) to
\[
(\mathbb{Z}^+ \sqcup \mathbb{R}^{++}, 1_{\mathbb{Z}^+}, \leq_{\Z})\,.
\]
\end{thms}

\begin{proof}
As observed in Corollary \ref{varecov}, the Cuntz semigroup of a
C$^*$-algebra $A$ of stable rank one is always the disjoint union
of the monoid $V(A)$ and $W(A)_+$.  Since $\Z$ is unital,
projectionless, and of stable rank one we have $V(\Z) \cong
\mathbb{Z}^+$.  By Proposition \ref{pureposcomparison} there is an
order-embedding
\[
\iota\colon W(\Z)_+ \to \mathbb{R}^{++}
\]
given by
\[
\iota(\langle a \rangle) = d_{\tau_{\Z}}(a),
\]
where $\tau_{\Z}$ is the unique normalised trace on $\Z$.  By
\cite[Theorem 2.1]{R1} there is a unital embedding of
$\mathrm{C}([0,1])$ into $\Z$ such that $\tau_{\Z}$ is implemented
by the uniform distribution on $[0,1]$.  Given $\lambda \in (0,1]$, let
$z_{\lambda} \in \mathrm{C}([0,1])$ be a positive function with
support $(0,\lambda)$.  It follows that
$d_{\tau_{\Z}}(z_{\lambda}) = \lambda$, whence $\iota$ is
surjective. We therefore have a bijection
\[
\varphi\colon W(\Z)=V(\Z)\sqcup W(\Z)_+\to \mathbb{Z}^+\sqcup \mathbb{R}^{++}\,.
\]
That $\varphi$ is an order-isomorphism follows from the fact that
$\Z$ has strict comparison of positive elements (\cite[Corollary 4.6]{R1}) and
Propositions \ref{projpureposcomp} and \ref{pureposcomparison}.
\end{proof}

\begin{nota}\label{zlam}
For each $\lambda \in (0,1]$ we will use $z_{\lambda}$ denote any
element in $\Z_{++}$ such that $d_{\tau_{\Z}}(z_{\lambda})=\lambda$.
\end{nota}

\begin{props}\label{iotadef}
Let $A$ be a C$^*$-algebra of stable rank one for which every trace is
faithful.  Then, the map
\[
\iota\colon W(A)_+\to \laff_b(\mathrm{T}(A))^{++}
\]
given by $\iota(\langle a\rangle)(\tau)=\dt(a)$ is a homomorphism.  If
$A$ has strict comparison of positive elements, then $\iota$ is an
order embedding.
\end{props}

\begin{proof}
The requirement that every trace on $A$ be faithful guarantees that
$\iota(\langle a \rangle)$ is strictly positive.
$A$ has stable rank one, so $W(A)_+$ is a semigroup by Proposition
\ref{pureposspec} and $\iota$ is a homomorphism.

If $A$ has strict comparison of positive elements, then $\iota$ is an
order embedding by Proposition \ref{pureposcomparison}.
\end{proof}

\begin{lms}\label{specialrep}
Let $A$ be a unital and $\mathcal{Z}$-stable C$^*$-algebra, with $a \in A_+$.
Then, $a$ is Cuntz equivalent to a positive element of the form $b \otimes \mathbf{1}_{\mathcal{Z}}
 \in A \otimes \mathcal{Z} \cong A$.
\end{lms}

\begin{proof}
Let $\psi: \mathcal{Z} \otimes \mathcal{Z} \to \mathcal{Z}$ be a $*$-isomorphism, and
put $\phi = (id_{\mathcal{Z}} \otimes \mathbf{1}_{\mathcal{Z}}) \circ \psi$.  By
\cite[Corollary 1.12]{TW1}, $\phi$ is approximately inner, and therefore so also is
\[
id_A \otimes \phi: A \otimes \mathcal{Z}^{\otimes 2} \to A \otimes \mathcal{Z} \otimes \mathbf{1}_{\mathcal{Z}}.
\]
In particular, there is a sequence of unitaries $u_n$ in $A \cong A \otimes \mathcal{Z}^{\otimes 2}$
such that 
\[
||u_n a u_n^* - \phi(a)|| \stackrel{n \to \infty}{\longrightarrow} 0. 
\]
Approximate unitary equivalence
preserves Cuntz equivalence classes, whence $\langle a \rangle = \langle \phi(a) \rangle$.
The image of $\phi(a)$ is, by construction, of the form $b \otimes \mathbf{1}_{\mathcal{Z}}$
for some $b \in A \otimes \mathcal{Z} \cong A$. 
\end{proof}

\begin{lms}\label{conelem}
Let $A$ be a unital, stably finite, and $\mathcal{Z}$-stable
$C^{*}$-algebra. Suppose that $f \in \laff (\mathrm{T}(A))^{++}$ is equal to $\dt(a)$ for some $a \in
\mathrm{M}_{\infty}(A)_{+}$. Then, the image of $a \otimes
z_{\lambda}$ in $\laff (\mathrm{T}(A \otimes \mathcal{Z}))^{++}$ is $\lambda \tilde{f}$, where
$\tilde{f} = \dt(a \otimes \mathbf{1}_{\mathcal{Z}})$.
\end{lms}

\begin{proof}
For any $\tau \in \mathrm{T}(A)$ one has
\begin{eqnarray*}
\dt(a \otimes z_{\lambda}) & = &\lim_{n \to \infty} \tau \left((a \otimes
z_{\lambda})^{1/n} \right)\\
& = & \lim_{n \to \infty} \tau(a^{1/n})\tau_{\Z}(z_{\lambda}^{1/n})\\
& = & \dt(a)d_{\tau_{\Z}}(z_{\lambda}) \\
& = & \lambda \dt(a).
\end{eqnarray*}

\end{proof}

\begin{cors}\label{cone}
Let $A$ be as in Lemma \ref{conelem}.  Then, the image of
$\mathrm{M}_{\infty}(A)_{+}$ under the map $\iota$ of Proposition \ref{iotadef}
 is a cone over $\mathbb{R}^{+}$
\end{cors}

\begin{proof}
It will be enough to prove that if $\lambda \in \mathbb{R}^+$ and $a \in A_+$, 
then there exists $b \in A_+$ with $\dt(b) = \lambda \dt(a)$.  Identify $A$ with
$A \otimes \mathcal{Z}$, and use Lemma \ref{specialrep} to find $b \in A_+$ such
that $b \otimes \mathbf{1}_{\mathcal{Z}}$ and $a$ are Cuntz equivalent.  It follows
that $\dt(b \otimes \mathbf{1}_{\mathcal{Z}}) = \dt(a)$ for each $\tau \in \mathrm{T}(A)$.  
Now $\dt(b \otimes z_{\lambda}) = \lambda \dt(a)$ by Lemma \ref{conelem}.
\end{proof}

Summarising, we have:

\begin{cors}
Let $A$ be a simple, unital, exact, finite, and $\Z$-stable C$^*$-algebra.
Then, the map $\iota$ of Proposition \ref{iotadef} is an order embedding,
and $W(A)_+$ is a $\mathbb{R}^+$-cone.
\end{cors}

\noindent
Note that exactness is required above in order to identify the image of $\iota$
with a collection of functions on $\mathrm{T}(A)$ as opposed to $\mathrm{QT}(A)$.

We close this section with an aside on some algebras of particular interest 
in Elliott's classification programme.
Recall that a C$^*$-algebra is said to have property (SP) if every hereditary
subalgebra contains a non-zero projection.  With Theorem \ref{WZ} in hand, we
can prove the following proposition:

\begin{props}\label{sp=smalltrace}
Let $A$ be a simple, unital, exact, finite, and $\mathcal{Z}$-stable C$^*$-algebra.
Then $A$ has property (SP) if and only if for every $\epsilon>0$ there exists
a non-zero projection $p \in A$ such that $\dt(p) = \tau(p) < \epsilon$ for every trace
on $A$.

In particular, a projection $p$ is Murray-von Neumann equivalent to a projection $q$
in a hereditary subalgebra $\overline{aAa}$ whenever $\tau \mapsto \dt(p)$ if uniformly
sufficiently small.
\end{props}

\begin{proof}
For the forward implication, write $A \cong A \otimes \mathcal{Z}$, and notice
that $\dt(1_A \otimes z_{\lambda}) = \lambda$, for all $\tau \in \mathrm{T}(A)$. Since $A$ has property
(SP), the algebra $\overline{(1_A \otimes z_{\lambda})A(1_A \otimes z_{\lambda})}$
contains a projection $p$, whence $p \precsim 1_A \otimes z_{\lambda}$.  Setting
$\lambda = \epsilon/2$, we have that
\[
\tau(p) \leq \dt(1_A \otimes z_{\epsilon/2}) < \epsilon, \ \text{for all } \tau \in \mathrm{T}(A).
\]

For the reverse implication, let $a \in A_+$ be given. The compactness of
$\mathrm{T}(A)$ and the lower semicontinuity of the function $f_a\colon \mathrm{T}(A) \to
\mathbb{R}^{++}$ given by $f_a(\tau) = \dt(a)$ (that follows from Proposition~\ref{lsc})
imply that there exists $\epsilon > 0$ such that $\dt(a) > \epsilon$, for every $\tau$ in $\mathrm{T}(A)$.
Choose a non-zero projection
$p$ in $A$ such that $\dt(p) = \tau(p) < \epsilon$ for every trace on $A$.  The hypotheses on
$A$ guarantee strict comparison for positive elements (cf.~\cite[Corollary 4.6]{R1}), so that
$p \precsim a$ inside $W(A)$.  Following the proof of Proposition~\ref{projcomplement},
we see that there is a projection $q \in \overline{aAa}$ which is Murray-von Neumann
equivalent to $p$.
\end{proof}

Let $\mathcal{B}$ be a class of unital C$^*$-algebras.  Recall that a unital C$^*$-algebra $A$
is said to be \emph{tracially approximately $\mathcal{B}$} (TA$\mathcal{B}$) if
for any $\epsilon > 0$, finite set $F \subset A$, and $a \in A_+$ there exists
a C$^*$-subalgebra $C$ of $A$ such that $C \in \mathcal{B}$, $\mathbf{1}_C \neq 0$,
and
\begin{enumerate}
\item $[f,\mathbf{1}_C] < \epsilon$, for all $f$ in $F$;
\item $\mathrm{dist}(\mathbf{1}_C f \mathbf{1}_C,C) < \epsilon$, for all $f$ in $F$;
\item $\mathbf{1}_A -\mathbf{1}_C$ is Murray-von Neumann equivalent to a projection in $\overline{aAa}$.
\end{enumerate}

One may wonder why the term ``tracially'' is used in the description of such algebras,
given that no reference to traces is made in their definition.  The reason is that
condition $(3)$ above can sometimes be replaced by the condition
\begin{enumerate}
\item[$(3)^{'}$]  $\tau(\mathbf{1}_C)> 1-\epsilon$, for all $\tau \in \mathrm{T}(A)$,
\end{enumerate}
provided that the class of TA$\mathcal{B}$ algebras is sufficiently well behaved.

TA$\mathcal{B}$ algebras are used mainly in Elliott's
classification program.  In this setting, it is necessary to
assume exactness, and the largest class for which classification
can be hoped for consists of $\mathcal{Z}$-stable algebras.  Since, in the
simple case, the program is more or less complete for purely
infinite algebras, we may also assume finiteness.  Taken together,
these conditions constitute the hypotheses of Proposition
\ref{sp=smalltrace}, and the proof of the proposition then shows
that conditions $(3)$ and $(3)^{'}$ above are equivalent.  Thus,
in most situations where TA$\mathcal{B}$ algebras might be useful,
there is no ambiguity in their definition.


\section{An embedding theorem}

In order to make sense of (EC) and (WEC), we must
define the categories in which the relevant invariants sit.

Let $\mathcal{I}$ denote the category whose objects are 4-tuples
\[
\left( (G_0,G_0^+,u),G_1,X,r \right),
\]
where $(G_0,G_0^+,u)$ is a simple
partially ordered Abelian group with distinguished order-unit
$u$ and state space $S(G_0,u)$, $G_1$ is a countable
Abelian group, $X$ is a metrizable Choquet simplex, and
$r\colon X \to S(G_0,u)$ is an affine map.  A morphism
\[
\Theta\colon \left( (G_0,G_0^+,u),G_1,X,r \right) \to
\left( (H_0,H_0^+,v),H_1,Y,s \right)
\]
in $\mathcal{I}$ is a 3-tuple
\[
\Theta = (\theta_{0}, \theta_{1}, \gamma)
\]
where
\[
\theta_{0}\colon (G_0,G_0^+,u) \to (H_0,H_0^+,v)
\]
is an order-unit-preserving positive homomorphism,
\[
\theta_{1}\colon G_1 \to H_1
\]
is any homomorphism, and
\[
\gamma\colon Y \to X
\]
is a continuous affine map that makes the diagram below commutative:
\[
\xymatrix{
{Y}\ar[r]^-{\gamma}\ar[d]^-{s} & {X}\ar[d]^-{r} \\
{S(H_0,v)}\ar[r]^-{\theta_0^*} & {S(G_0,u)\,.}}
\]
For a simple unital $C^{*}$-algebra $A$ the Elliott invariant
$I(A)$ is an element of $\mathcal{I}$, where
$(G_0,G_0^+,u) = (\mathrm{K}_{0}(A), \mathrm{K}_{0}(A)^{+}, [1_{A}])$,
$G_1  =  \mathrm{K}_{1}(A)$, $X  =  \mathrm{T}(A)$,
and $r_A$ is given by evaluating a given trace at a $\mathrm{K}_0$-class.
Given a class $\mathcal{C}$ of simple unital C$^*$-algebras, let
$\mathcal{I}(\mathcal{C})$ denote the subcategory of $\mathcal{I}$ whose
objects can be realised as the Elliott invariant of a member of
$\mathcal{C}$, and whose morphisms
are all admissible maps between the now specified objects.

The definition of $\mathcal{I}$ removes an ambiguity from
the statement of (EC), namely, what is meant
by an isomorphism of Elliott invariants.  We now do the same for (WEC).
Let $\mathcal{W}$ be the category whose objects are ordered pairs
\[
\left( (W(A),\langle 1_A \rangle), I(A) \right),
\]
where $A$ is a simple, unital, exact, and stably finite C$^*$-algebra, $(W(A),\langle 1_A \rangle)$
is the Cuntz semigroup of $A$ together with the distinguished order-unit $\langle 1_A \rangle$,
and $I(A)$ is the Elliott invariant of $A$. A morphism
\[
\Psi\colon \left( (W(A),\langle 1_A \rangle), I(A) \right) \to
\left( (W(B),\langle 1_B \rangle), I(B) \right)
\]
in $\mathcal{W}$ is an ordered pair
\[
\Psi = (\Lambda, \Theta),
\]
where $\Theta = (\theta_0,\theta_1,\gamma)$ is a morphism  in $\mathcal{I}$ and
$\Lambda\colon (W(A),\langle 1_A\rangle)\to (W(B),\langle 1_B\rangle)$ is an order- and
order-unit-preserving semigroup homomorphism satisfying two
compatibility conditions:  first,
\[
\xymatrix{
 { (V(A),\langle 1_A \rangle)}\ar[r]^-{\Lambda|_{V(A)}}\ar[d]^-{\rho}&
{(V(B),\langle 1_B \rangle)}\ar[d]^-{\rho} \\
(\mathrm{K}_0(A),[1_A])\ar[r]^-{\theta_0}&(\mathrm{K}_0(B),[1_B])\,,}
\]
where $\rho$ is the usual Grothendieck map from $V(\bullet)$ to $\mathrm{K}_0(\bullet)$ (recall
that there is an order-unit-preserving order-embedding of $(V(A),\langle 1_A \rangle)$
into $(W(A),\langle 1_A \rangle)$, and that Cuntz equivalence of projections agrees with Murray-von Neumann
equivalence in stably finite algebras);
second,
\[
\xymatrix{
\mathrm{LDF}(B)\ar[r]^-{\Lambda^*}\ar[d]^-{\eta} & \mathrm{LDF}(A)\ar[d]^-{\eta} \\
\mathrm{T}(B)\ar[r]^-{\gamma}&\mathrm{T}(A)\,,}
\]
where $\eta$ is the affine bijection between
$\mathrm{LDF}(\bullet)$ and $\mathrm{T}(\bullet)$ given by
$\eta(\dt) = \tau$ (see~\cite[Theorem II.2.2]{bh}).  These
compatibility are automatically satisfied if $\Psi$ is induced by
a $*$-homomorphism $\psi: A \to B$.

Recall that we have previously defined, for a C$^*$-algebra with stable rank one, a semigroup homomorphism
\[
\iota\colon W(A)_+ \to \laff_b(\mathrm{T}(A))^{++}
\]
by
\[
\iota(\langle a \rangle)(\tau)=\dt(a), \ \text{for all } \tau \in \mathrm{T}(A)\,.
\]
In the following definition we generalise the semigroup and order
structure on \mbox{$\mathbb{Z}^+\sqcup\mathbb{R}^{++}$} considered
in Theorem \ref{WZ}. Semigroups of this type have
been considered previously in the study of multiplier algebras
(see~\cite{percan}).

\begin{dfs}\label{wtilde}
Let $A$ be a unital C$^*$-algebra. Define a semigroup structure on the set
\[
\wt:= V(A) \sqcup \laff_b(\mathrm{T}(A))^{++}
\]
by extending the natural semigroup operations and setting $[p]+f=\widehat{p}+f$, where $\widehat{p}(\tau)=\tau(p)$.
Define an order $\leq$ on $\wt$ such that:
\begin{enumerate}[{\rm (i)}]
\item $\leq$ agrees with the usual order on $V(A)\,;$
\item $f \leq g$ for $f$, $g$ in $\laff(\mathrm{T}(A))^{++}$ if and only if
\[
f(\tau) \leq_{\mathbb{R}} g(\tau) \ \text{for all } \tau \in \mathrm{T}(A)\,;
\]
\item $f \leq [p]$ for
$[p] \in V(A)$ and $f$ in $\laff (\mathrm{T}(A))^{++}$ if and only if
\[
f(\tau) \leq_{\mathbb{R}} \tau(p) \ \text{for all } \tau \in \mathrm{T}(A)\,;
\]
\item $[p] \leq f$ for $f$, $[p]$ as in $(iii)$ whenever
\[
\tau(p) <_{\mathbb{R}} f(\tau) \ \text{for all } \tau \in \mathrm{T}(A)\,.
\]
\end{enumerate}

\end{dfs}

Let $\mathcal{\widetilde{W}}$ be the category whose objects are
of the form $(\widetilde{W}(A),[1_A])$ for some exact, unital, and stable rank
one C$^*$-algebra $A$, and whose morphisms are positive order-unit-preserving
homomorphisms
\[
\Gamma\colon (\widetilde{W}(A),[1_A]) \to (\widetilde{W}(B),[1_B])
\]
such that
\[
\Gamma(V(A)) \subseteq V(B)
\]
and
\[
\Gamma|_{\laff_b(\mathrm{T}(A))^{++}}\colon \laff_b(\mathrm{T}(A))^{++} \to \laff_b(\mathrm{T}(B))^{++}
\]
is induced by a continuous affine map from $\mathrm{T}(B)$ to $\mathrm{T}(A)$.

For the next definition, we remind the reader that $V(A) \cong\mathrm{K}_0(A)^+$
for a C$^*$-algebra of stable rank one.

\begin{dfs}\label{recovfunc} Let $\mathcal{C}$ denote the class of simple, unital, exact, and
stable rank one C$^*$-algebras.
Let
\[
F\colon\mathbf{Obj}(\mathcal{I}(\mathcal{C})) \to \mathbf{Obj}(\mathcal{\widetilde{W}})
\]
be given by
\[
F\left( (\mathrm{K}_0(A),\mathrm{K}_0(A)^+,[1_A]),\mathrm{K}_1(A),\mathrm{T}(A),r_A \right) =
(\widetilde{W}(A),[1_A]).
\]
Define
\[
F\colon\mathbf{Mor}(\mathcal{I}(\mathcal{C})) \to \mathbf{Mor}(\mathcal{\widetilde{W}})
\]
by sending $\Theta = (\theta_0,\theta_1,\gamma)$ to the morphism
\[
\Gamma\colon(\widetilde{W}(A),[1_A]) \to (\widetilde{W}(B),[1_B])
\]
given by $\theta_0$ on $\mathrm{K}_0(A)^+ = V(A)$ and induced by $\gamma$ on $\laff_b(\mathrm{T}(A))^{++}$.
\end{dfs}

The next proposition holds by definition.

\begin{props} With $\mathcal{C}$ as in Definition \ref{recovfunc}, the map
 $F\colon\mathcal{I}(\mathcal{C}) \to \mathcal{\widetilde{W}}$ is a functor.
\end{props}

For the theorem below, we remind the reader that the definition of the map 
$\iota$ is contained in Proposition \ref{iotadef}.

\begin{thms}\label{embedding}
Let $A$ be a simple, unital, and exact C$^*$-algebra having stable rank one and strict comparison
of positive elements.  Then, there is an order embedding
\[
\phi\colon W(A) \to \wt
\]
such that $\phi|_{V(A)} = \mathrm{id}_{V(A)}$ and $\phi|_{W(A)_+} = \iota$.
\end{thms}

\begin{proof}
The map $\phi$ is well-defined, so it will suffice to prove that
it is an order embedding.  We verify conditions (i)-(iv) from
Definition \ref{wtilde}:
the image of $\phi|_{V(A)}$ is $V(A)$, with the same order, so (i)
is satisfied;  (ii) and (iii) follow from Proposition~\ref{pureposcomparison};
(iv) is Proposition~\ref{projpureposcomp}.
\end{proof}

We are now ready to prove Theorem \ref{ecwec1}. 

\begin{thms}
(EC) implies (WEC) for the class of simple, unital, separable, and nuclear
C$^*$-algebras which absorb the Jiang-Su algebra $\mathcal{Z}$ tensorially.
\end{thms}

\begin{proof}
Algebras in the class under consideration are either purely infinite
or stably finite (cf. \cite{gjsu}).
The theorem is trivial for the subclass of purely infinite algebras,
due to the degenerate nature of the Cuntz semigroup in this setting.
The remaining case is that of stable rank one.

Let $A$ and $B$ be simple, separable, unital, nuclear, and stably finite
C$^*$-algebras with strict comparison of positive elements, and suppose
that (EC) holds.  Let there be given an isomorphism
\[
\phi\colon \left(W(A), \langle 1_A \rangle, I(A)\right) \to
\left(W(B), \langle 1_B \rangle, I(B)\right).
\]
Then by restricting $\phi$ we have an isomorphism
\[
\phi|_{I(A)}\colon I(A) \to I(B),
\]
and we may conclude by (EC) that there is a $*$-isomorphism $\Phi\colon A \to B$
such that $I(\Phi) = \phi|_{I(A)}$.  $\Phi$ is unital and so preserves the Cuntz class
of the unit.  The compatibility conditions imposed on $\phi$ (see the discussion preceding 
Definition \ref{wtilde}) together with Theorem
\ref{embedding} ensure that $\phi|_{W(A)}$
is determined by $\phi|_{V(A)}$ and $\phi^{\sharp}: \mathrm{T}(B) \to \mathrm{T}(A)$.
Thus, $\Phi$ induces $\phi$, and (WEC) holds.
\end{proof}

\noindent
Note that the semigroup homomorphism $\phi$ in Theorem~\ref{embedding} is an isomorphism if and only if
$\iota$ is surjective.

Let $\mathrm{(EC)}^{'}$ and $\mathrm{(WEC)}^{'}$ denote the statements (EC) and (WEC),
respectively, but expanded to apply to all simple, unital, exact, and stably finite C$^*$-algebras.
Collecting the results of this section we have:

\begin{thms}\label{ecwec}
Let $\mathcal{C}$ be a class of simple, unital, exact, finite, and $\mathcal{Z}$-stable
C$^*$-algebras.  Suppose that $\iota$ is surjective for each member of $\mathcal{C}$.
Then, $(EC)^{'}$ and $(WEC)^{'}$ are equivalent in $\mathcal{C}$.  Moreover, there is a functor
$G\colon \mathcal{I}(\mathcal{C}) \to \mathcal{W}$ such that
\[
G(I(A)) \stackrel{\mathrm{def}}{=} \left( F(I(A)),I(A) \right) = 
((\widetilde{W}(A),[1_A]),I(A)) \cong ((W(A),\langle 1_A \rangle),I(A)).
\]
\end{thms}

Even in situations where (EC) holds, there is no 
inverse functor which reconstructs C$^*$-algebras from Elliott invariants.  (This
is {\it not} the same as saying that one cannot reconstruct the algebra from
the Elliott invariant at all --- this is always possible when one has a range
result for a class of algebras satisfying (EC).)
Contrast this with Theorem \ref{ecwec}, where $G$ reconstructs the finer invariant
from the coarser one functorially.  

We now see that (EC) and (WEC) are equivalent
among simple, unital, separable, nuclear, finite, and $\mathcal{Z}$-stable C$^*$-algebras 
whenever $\iota$ is surjective.  (It is not clear whether the converse holds.)
In Sections 6, 7, and 8 we will prove that $\iota$ is surjective for algebras
satisfying hypotheses (i), (ii), or (iii) of Theorem \ref{equiv}, respectively, 
thereby proving the theorem.

We note that if $\iota$ is surjective and $A$ satisfies the hypotheses of Theorem \ref{ecwec},
then the invariant
\[
((W(A),\langle 1_A \rangle),I(A))
\]
carries redundant information.  $A$ has stable rank one, so one
may, by using Corollary \ref{varecov}, recover $V(A)
\cong\mathrm{K}_0(A)^+$, and hence
$(\mathrm{K}_0(A),\mathrm{K}_0(A)^+,[1_A])$, from $(W(A),\langle
1_A \rangle)$.  The convex affine space $\mathrm{T}(A)$ is identified
with $\mathrm{LDF}(A)$ (although we cannot, in general, recover
the topology on $\mathrm{T}(A)$ -- see the discussion following
Corollary I.2.2 of \cite{bh}). The pairing
$r_A$ can be recovered by applying the elements of
$\mathrm{LDF}(A)$ to $V(A)\cong\mathrm{K}_0(A)^+$.

We close this section by observing that if $\iota$ is surjective, then the failure of
the order on $W(\bullet)$ to be algebraic in general is easily
explained.

\begin{props}
Let $A$ be an exact C$^*$-algebra with strict comparison of positive elements.  Suppose
that $\iota$ is surjective and that each $\tau \in \mathrm{T}(A)$ is faithful.
Let $a \precsim b$ in $\mathrm{M}_{\infty}(A)_{++}$.  Then,
there exists a positive element $c \in \mathrm{M}_{\infty}(A)_{++}$ such that
$a \oplus c \sim b$ if and only if the difference
\[
\dt(b)-\dt(a)\colon \mathrm{T}(A) \to \mathbb{R}^+
\]
is in $\laff_b(\mathrm{T}(A))^{++}$.
\end{props}

\begin{proof}
If $b \sim a \oplus c$, then $\dt(b) - \dt(a) = \dt(c)$ and $\dt(c) \in \laff_b(\mathrm{T}(A))^{++}$ by
Proposition \ref{lsc}.

Suppose that $f(\tau) := \dt(b)-\dt(a) \in \laff_b(\mathrm{T}(A))^{++}$.  Choose, by the surjectivity of
$\iota$, an element $c \in \mathrm{M}_{\infty}(A)_{++}$ for which $\dt(c) = f(\tau)$.
Then $\dt(a \oplus c) = \dt(b)$, whence $a \oplus c \sim b$ by Proposition
\ref{pureposcomparison}.
\end{proof}


\section{The structure of $\mathrm{K}_0^*$}

The Grothendieck enveloping group of
$W(A)$ is denoted $\mathrm{K}_0^*(A)$, and its structure has been
previously analysed in \cite{bh}, \cite{Cu}, \cite{han},
 and~\cite{P1}.
Because $W(A)$ carries its own order coming from the Cuntz
comparison relation, $\mathrm{K}_0^*(A)$ may be given two natural
(partial) orderings. For an abelian semigroup $M$ with a partial
order $\leq$ that extends the algebraic order, we use $G(M)$ to
denote its enveloping group. Write $\gamma\colon M\to G(M)$ for
the natural Grothendieck map. We define the following cones:
\[
G(M)^+=\gamma (M)\,,
\]
and
\[
G(M)^{++}=\{\gamma(x)-\gamma(y)\mid x,\,y\in M\text{ and }y\leq
x\}\,.
\]
Since $M$ is partially ordered, so is $(G(M), G(M)^{++})$.
Clearly, $G(M)^+\subseteq G(M)^{++}$, and the inclusion may be
strict. Therefore, $(G(M),G(M)^+)$ is also partially ordered. For
the reader's convenience, we offer a short argument which shows
the cone $G(M)^{++}$ to be strict (compare with~\cite{han}
and~\cite{bh}). Assume that $\gamma(x)-\gamma(y)\in
G(M)^{++}\cap(-G(M)^{++})$. Then there are elements $s$, $t$, $u$,
$v$ in $M$ such that
\[
x+z\leq y+z\,,\,\, t+v\leq s+v\,,\,\, x+s+u=y+t+u\,,
\]
so that $\gamma(y)-\gamma(x)=\gamma(s)-\gamma(t)\in G(M)^{++}$. Set $w=u+v+z+t$
and check that $x+w=y+w$, whence $\gamma(x)=\gamma(y)$.

Recall that a partially ordered Abelian group with order-unit
$(G,G^+,u)$ is \emph{Ar\-chi\-me\-dean} provided that $nx\leq y$ for $x$, $y\in G$ and for all natural numbers
$n$ only if $x=0$ (see~\cite[p. 20]{G}). This is
equivalent (by~\cite[Theorem 4.14]{G}) to saying that the order on
$G$ is determined by its states, i.e., $G^+=\{x\in G\mid
s(x)\geq 0\text{ for all }s\in \mathrm{S}(G,u)\}$. (Recall that a state
$s$ on $(G,G^+,u)$ is a positive group homomorphism into $\mathbb{R}$
such that $s(u)=1$ --- $s$ need not be order preserving, in contrast with
a state on a positive ordered Abelian semigroup.) We say that
$(G,G^+)$ is \emph{unperforated} if $nx\geq 0$ implies that $x\geq
0$ (see~\cite{G}). Archimedean directed groups are
unperforated (cf.~\cite[Proposition 1.24]{G}).

For an element $a$ in $M_{\infty}(A)_+$, we shall denote by $[a]$
the class of $\langle a\rangle$ in $\mathrm{K}_0^*(A)$.

\begin{lms}\label{k*comparison}
Let $A$ be a simple C$^*$-algebra with strict comparison of positive
elements. Suppose that $\mathrm{M}_{\infty}(A)_{++} \neq \emptyset$.
Then:
\[
\mathrm{K}_0^*(A)^{++}=\{[a]-[b]\mid a,\, b\in
M_{\infty}(A)_+\text{ and }\dt(a)\geq\dt (b)\text{ for all
}\tau\in\mathrm{QT}(A)\}\,.
\]
\end{lms}

\begin{proof}
By the properties of dimension functions, it is clear that if $a$,
$b\in M_{\infty}(A)_+$ and $b\precsim a$, we have
$\dt(b)\leq\dt(a)$ for any $\tau\in \mathrm{QT}(A)$.

For the converse inclusion, let $[a]-[b]\in \mathrm{K}_0^*(A)$ be
such that $\dt(b)\leq \dt(a)$ for each $\tau\in\mathrm{QT}(A)$.
Then, for any $0 \neq c \in \mathrm{M}_{\infty}(A)_{++}$ we have
$a \oplus c, b \oplus c \in \mathrm{M}_{\infty}(A)_{++}$ and
\[
\dt(b \oplus c) \leq \dt(a \oplus c)\,.
\]
It follows from Proposition~\ref{pureposcomparison} that
\[
b \oplus c \precsim a \oplus c\,,
\]
and thus $[a]-[b]=[a\oplus c]-[b\oplus c]\in
\mathrm{K}_0^*(A)^{++}$.
\end{proof}

\begin{cors}\label{k*archimedean}
Let $A$ be a C$^*$-algebra satisfying the hypotheses of
Lemma~\ref{k*comparison}. Then
$(\mathrm{K}_0^*(A),\mathrm{K}_0^*(A)^{++})$ is Archimedean, and
in particular is unperforated.
\end{cors}

\begin{proof}
The second conclusion follows from the first since, as observed
above, ar\-chi\-me\-dean groups are unperforated. (Notice that
$\mathrm{K}_0^*(A)$ is directed since $A$ is unital.)

We only need to show that if $[a]-[b]\in \mathrm{K}_0^*(A)$ is
such that $s([a]-[b])\geq 0$ for any state $s$ on
$\mathrm{K}_0^*(A)$ (i.e. $s([b])\leq s([a])$), then $[a]-[b]\in
\mathrm{K}_0^*(A)^{++}$. Recalling that the states on
$\mathrm{K}_0^*(A)$ are precisely the dimension functions, we have
that in particular $\dt(b)\leq\dt(a)$ for any quasitrace $\tau$,
hence we may use Lemma~\ref{k*comparison}.
\end{proof}

We shall show below that $\mathrm{K}_0^*(A)$ is also unperforated
when endowed with the ordering defined by taking as positive cone
$\mathrm{K}_0^*(A)^+=\gamma (W(A))$, that is, the image of $W(A)$
under the Grothendieck map.

A partially ordered semigroup $(M,\leq)$ is said to be
\emph{almost unperforated} if for all $x$, $y$ in $M$ and
$n\in\mathbb{N}$ with $(n+1)x\leq ny$, one has that $x\leq y$. A
simple partially ordered group $(G,G^+)$ is \emph{weakly
unperforated} if $nx\in G^+\setminus\{0\}$ implies that $x\in
G^+\setminus\{0\}$ (\cite[Lemma 3.4]{R1}).

\begin{props}\label{k*unperf}
Let $A$ be a simple, unital, exact, and finite C$^*$-algebra which absorbs
the Jiang-Su algebra $\mathcal{Z}$ tensorially. Then, the
partially ordered Abelian group
$(\mathrm{K}_0^*(A),\mathrm{K}_0^*(A)^+)$ is weakly
unperforated.
\end{props}

\begin{proof}
We have already noticed that $A$ has strict comparison of positive
elements, by Corollary 4.6 of~\cite{R1}. The simplicity of $A$
guarantees that each trace on $A$ is faithful. Since $1_A \otimes
z_1 \in A \otimes \mathcal{Z} \cong A$, we have that
$\mathrm{M}_{\infty}(A)_{++} \neq \emptyset$.  Thus, $A$ satisfies
the hypotheses of Lemma~\ref{k*comparison}.

Given $[a] \in \mathrm{K}_0^*(A)^+$, for $a\in M_{\infty}(A)_+$,
we may assume that $a \in \mathrm{M}_{\infty}(A)_{++}$. 
To see this, first identify $A$ with $A \otimes \mathcal{Z}$, and
replace $a$ with a Cuntz equivalent element $b \otimes \mathbf{1}_{\mathcal{Z}}$
(see Lemma \ref{specialrep}).
Now for each $\tau \in \mathrm{T}(A)$ we have
\[
\dt(a) = \dt(b \otimes \mathbf{1}_{\mathcal{Z}}) = \dt(b \otimes z_1)
\]
(see Notation \ref{zlam} and Lemma \ref{conelem}).
Now $[a] = [b \otimes z_1]$ by
Lemma~\ref{k*comparison} and the proof of the fact that
$\mathrm{K}_0^*(A)^{++}$ is strict.  We have $z_1 \in \mathcal{Z}_{++}$
by construction, and a straighforward functional calculus argument then
shows that $b \otimes z_1 \in \mathrm{M}_{\infty}(A)_{++}$.

Suppose that $[a]$, $[b]\in\mathrm{K}_0^*(A)^+$ are such that
\[
(n+1)[a]\leq n[b]\,, \ \text{for some } n \in \mathbb{N}\,.
\]
This means that there is $c\in M_{\infty}(A)_+$ such that $(n+1)[a]+[c]=n[b]$.

Assume that $a,b \in \mathrm{M}_{\infty}(A)_{++}$. By
Lemma~\ref{k*comparison}, we have $(n+1)\dt(a)+\dt(c)= n\dt(b)$,
whence $\dt(a)+\frac{1}{n}\dt(a\oplus c) =\dt(b)$. Invoke
Corollary \ref{cone} to find a (purely positive) element $c'$
such that $\frac{1}{n}\dt(a\oplus c)=\dt(c')$. Now,
Proposition~\ref{pureposcomparison} implies that $a \oplus c' \sim
b$, whence $[a]+[c']= [b]$.  This shows that $\mathrm{K}_0^*(A)^+$
is almost unperforated. Apply Lemma 3.4 of~\cite{R1} and the
discussion thereafter to conclude that
$(\mathrm{K}_0^*(A),\mathrm{K}_0^*(A)^+)$ is weakly unperforated.
\end{proof}

Note that if $A$ is simple, then
$(\mathrm{K}_0^*(A),\mathrm{K}_0^*(A)^+)$ is a simple group. This
raises the question of whether
$(\mathrm{K}_0^*(A),\mathrm{K}_0^*(A)^{++})$ will also be simple
for a simple C$^*$-algebra $A$. We give a criterion below to
decide when a given (positive) element in $\mathrm{K}_0^*(A)^{++}$
is an order-unit. If $a\in M_{\infty}(A)_{+}$, write $n\cdot a$ to mean $a\oplus\cdots\oplus a$ ($n$ times).

\begin{props}
Let $A$ be a unital, simple, stably finite, exact C$^*$-algebra with strict comparison of positive elements. Suppose that
$M_{\infty}(A)_{++}$ is non-empty. Then, an element $[a]-[b]\in
\mathrm{K}_0^*(A)^{++}$ is an order-unit if and only if there is
$\epsilon>0$ such that $\dt(a)-\dt(b)>\epsilon$ for all traces
$\tau$.
\end{props}

\begin{proof}
If $[a]-[b]$ is an order-unit, then clearly $[a]\neq 0$. 
If $b=0$ then 
\[
\dt(a)-\dt(b)=\dt(a)>0
\]
for each $\tau \in \mathrm{T}(A)$.
The function $\tau \mapsto \dt(a)$ is lower semicontinuous on a compact
set, and therefore achieves a minimum $\delta > 0$.  Setting $\epsilon=\delta/2$
gives the desired conclusion.  

Now suppose that $b \neq 0$.
There is a natural number $n$ such that $[a]\leq n[a]-n[b]$, hence
we can find $c\in M_{\infty}(A)_+$ such that $a\oplus c\oplus
n\cdot b\precsim n\cdot a\oplus c$.
Therefore, for any $\tau\in\mathrm{T}(A)$, we have
$\dt(a)+n\dt(b)\leq n\dt(a)$.  Since $b \neq 0$ we conclude that
\[
(n-1)(\dt(a)-\dt(b))>\dt(b)>0.
\]
Using the same argument as in the $b=0$ case, we conclude that there
is some $\epsilon>0$ such that $\dt(b) > (n-1)\epsilon$ for every $\tau
\in \mathrm{T}(A)$.  It follows that $\dt(a)-\dt(b) > \epsilon$, as desired.

Conversely, if $\dt(a)-\dt(b)>\epsilon$ for all $\tau$, choose $n$
such that $\dt(n\cdot a)-\dt(n\cdot b)=n(\dt(a)-\dt(b))>1=\dt(1_A)$. Let
$c\in M_{\infty}(A)_{++}$. Then
\[
\dt(n\cdot a\oplus n\cdot c)-\dt(n\cdot b\oplus n\cdot c)>\dt(1_A)\,,
\]
whence $\dt(n\cdot a\oplus n\cdot c)>\dt(n\cdot b\oplus n\cdot c\oplus 1_A)$ for all
$\tau$. If follows now from Proposition~\ref{pureposcomparison} that
$n\cdot b\oplus n\cdot c\oplus 1_A\precsim n\cdot a\oplus n\cdot c$. This implies that
$n([a]-[b])\geq [1_A]$, whence $[a]-[b]$ is an order-unit.
\end{proof}

\begin{lms}
\label{k*struc} Let $A$ be a C$^*$-algebra with stable rank one
and such that the semigroup $W(A)_+$ of purely positive elements
is non-empty. Then there exists an ordered group isomorphism
\[
\alpha\colon (\mathrm{K}_0^*(A), \mathrm{K}_0^*(A)^{++})\to
(G(W(A)_+),G(W(A)_+)^+)\,.
\]
If, furthermore, $A$ is simple and $\Z$-stable, then $\alpha([1_A])=([1\otimes z_1])$.
\end{lms}

\begin{proof}
Recall from Section 2 that if $A$ has stable rank one, then
$W(A)=V(A)\sqcup W(A)_+$. Denote by $\gamma\colon W(A)_+\to
G(W(A)_+)$ the Grothendieck map, and choose any element $c\in
W(A)_+$. Then, define
\[
\alpha\colon W(A)\to G(W(A)_+)
\]
by $\alpha(\langle a\rangle)=\gamma(\langle a\rangle)$ if $\langle
a\rangle\in W(A)_+$, and by $\alpha(\langle
p\rangle)=\gamma(\langle p\rangle+c)-\gamma(c)$ for any projection
in $M_{\infty}(A)$.

Note that $\alpha$ is a well defined semigroup homomorphism.
Indeed, since $A$ has stable rank one, $\langle p\rangle +c\in
W(A)_+$ whenever $c\in W(A)_+$ (Lemma \ref{varecov}), and if $c'\in W(A)_+$ is any other
element, then one has that $\gamma(\langle
p\rangle+c)-\gamma(c)=\gamma(\langle p\rangle+c')-\gamma(c')$.

In order to check that $\alpha$ is a homomorphism, let $p$, $q$
and $a$ be elements in $M_{\infty}(A)_+$ with $p$ and $q$
projections and $a$ purely positive. Then,
\begin{eqnarray*}
\alpha(\langle
p\rangle+\langle q\rangle) & = & \gamma(\langle p\oplus
q\rangle+2c)-\gamma(2c) \\
& = & \gamma(\langle
p\rangle+c)-\gamma(c)+\gamma(\langle
q\rangle+c)-\gamma(c) \\
& = & \alpha(\langle p\rangle)+\alpha(\langle
q\rangle).
\end{eqnarray*}

\noindent
Also
\begin{eqnarray*}
\alpha(\langle p\rangle+\langle a\rangle)& = & \gamma(\langle p\oplus a\rangle) \\
& = & \gamma(\langle p\oplus a\rangle+c)-\gamma(c) \\
& = & \gamma(\langle p\rangle+c)-\gamma(c)+\gamma(\langle a\rangle) \\
& = & \alpha(\langle p\rangle)+\alpha(\langle a\rangle).
\end{eqnarray*}

By definition,
\[
G(W(A)_+)^+ = \{ [x]-[y] \ | \ x,y \in W(A)_+, \ y +r \leq x + r \ \mathrm{for} \ \mathrm{some} \ r \in W(A)_+ \},
\]
whence $\alpha(W(A)_+) = \gamma(W(A)_+) \subseteq G(W(A)_+)^+$ by construction.  If
$p \in \mathrm{M}_{\infty}(A)$ is a projection, then its image under alpha
is $\gamma(\langle p \rangle +c) - \gamma(c)$.  Since $c \leq \langle p \rangle + c$,
we conclude that $\alpha(\langle p \rangle) \in G(W(A)_+)^+$, too. 
Thus, $\alpha (W(A))\subseteq G(W(A)_+)^+$, 
and so $\alpha$ extends to an ordered group homomorphism
\[
\alpha\colon \mathrm{K}_0^*(A)=G(W(A))\to G(W(A)_+)\,,
\]
given by the rule $\alpha([a]-[b])=\alpha(\langle
a\rangle)-\alpha(\langle b\rangle)$. Evidently, $\alpha$ is
surjective and satisfies
\[\
\alpha (\mathrm{K}_0^*(A)^{++})\subseteq G(W(A)_+)^+\,
\]

To prove injectivity, assume that $\alpha (\langle
a\rangle)=\alpha(\langle p\rangle)$ for $\langle a\rangle\in
W(A)_+$ and $p$ a projection. This means that $\gamma(\langle
a\rangle)=\gamma(\langle p\rangle+c)-\gamma(c)$, and hence
$\langle a\rangle+c+c'=\langle p\rangle+c+c'$ for some $c'\in
W(A)$. Thus $[a]=[p]$ in $\mathrm{K}_0^*(A)$. If, for projections
$p$ and $q$, we have that $\alpha(\langle p\rangle)=\alpha(\langle
q\rangle)$, then $\gamma(\langle
p\rangle+c)-\gamma(c)=\gamma(\langle q\rangle+c)-\gamma(c)$.  It
follows that $[p]=[q]$ in $\mathrm{K}_0^*(A)$.

Finally, if $A$ is simple and $\Z$-stable, we may apply Proposition~\ref{pureposcomparison} to conclude that
\[
(1_A\otimes 1_{\mathcal{Z}})\oplus (1_A\otimes z_1)\sim (1_A\otimes
z_1)\oplus (1_A\otimes z_1).
\]
Thus, there is an identification of $A$ with $A \otimes \mathcal{Z}$ for which
\[
\alpha([1_A])=\gamma (\langle
(1_A\otimes 1_{\mathcal{Z}})\oplus (1_A\otimes z_1)\rangle)-\gamma(\langle
1_A\otimes z_1\rangle)=\gamma (\langle 1_A\otimes
z_1\rangle)=\alpha ([1_A\otimes z_1]).
\]
(Note that the $1_A$s on the far right and far left are, strictly speaking, not the
same.)
\end{proof}

\begin{cors}\label{k*structure}
Let $A$ be simple, unital, and exact C$^*$-algebra having stable rank one
and strict comparison of positive elements.  Suppose further that $\mathrm{M}_{\infty}(A)_{++} \neq \emptyset$.
Then, $\mathrm{K}_0^*(A)$ is the Grothendieck
enveloping group of $\iota(W(A)_+)$, where $\iota$ is the map
defined in Proposition~\ref{iotadef}.
\end{cors}

\begin{proof}
Under the hypotheses, $\iota$ is an order-embedding (see
Theorem~\ref{embedding}). The result then follows from
Lemma~\ref{k*struc}.
\end{proof}

\noindent Corollary \ref{k*structure} gives a version of Theorem
III.3.2 of \cite{bh} for C$^*$-algebras which may lack non-trivial
projections.

We close this section summarizing our findings in the following:

\begin{thms}
Let $A$ be a simple, unital, nuclear and finite C$^*$-algebra which is $\mathcal{Z}$-stable. Then,
\begin{enumerate}[{\rm (i)}]

\item $(\mathrm{K}_0^*(A),\mathrm{K}_0^*(A)^{++})$ is an Archimedean partially ordered Abelian group.

\item $(\mathrm{K}_0^*(A),\mathrm{K}_0^*(A)^+)$ is a simple and weakly unperforated partially ordered Abelian group.

\item $\mathrm{K}_0^*(A)=G(\iota(W(A)_+)$, where $\iota\colon W(A)_+\to\laff_b(\mathrm{T}(A))^{++}$ is defined as in~\ref{iotadef}.
\end{enumerate}
\end{thms}

\section{$\Z$-stable algebras with finitely many pure tracial states}

In the final sections of the paper, we study the surjectivity of
the order-embedding $\iota$. In  this section we study algebras which
satisfy the hypotheses of Theorem \ref{equiv} by way of having finitely many pure
tracial states.
We begin by establishing a closure property for the image of $\iota$.

\begin{lms}
\label{specrep} Let $A$ be a simple, unital, exact, finite,
and $\Z$-stable C$^*$-algebra:
\[
A \stackrel{\phi}{\cong} A \otimes \mathcal{Z},  
\]
where $\phi$ is as in the proof of Lemma \ref{specialrep}.
Suppose that $a \in
\mathrm{M}_{\infty}(A)_{+}$ is such that $\dt(a) \leq r$, for some
$r \in \mathbb{R}^{++}$ and for all $\tau\in\mathrm{T}(A)$. Then,
for any $z$ in $\Z$ such that $z\sim z_r$, there exists $\tilde{a}
\in \mathrm{M}_{\infty}(A)_+$ such that
\[
a \sim \tilde{a} \leq (1_A \oplus 1_A) \otimes z \in \mathrm{M}_{\infty}(A \otimes \Z)_+ 
\stackrel{\phi}{\equiv} \mathrm{M}_{\infty}(A)_+.
\]
\end{lms}

\begin{proof}
We assume throughout the proof that whenever elementary tensors in $A \otimes \Z$
are mentioned, they are being identified with elements of $A$ via $\phi$.

Suppose first that $a \sim p$ for some projection $p \in \mathrm{M}_{\infty}(A)$.
Since
\[
\dt(a) \leq r < 2r = \dt((1_A \oplus 1_A) \otimes z)\,, \
\text{for all } \tau \in \mathrm{T}(A)\,,
\]
we have that $a \sim p \precsim (1_A \oplus 1_A) \otimes z$ by
Proposition \ref{projpureposcomp}. Applying~\cite[Proposition
2.4]{Rfunct} we may find $x \in \mathrm{M}_{\infty}(A)$ such that
\[
x^*\left((1_A \oplus 1_A) \otimes z\right)x = (p-\epsilon)_+ \sim p
\sim a,
\]
so that $\tilde{a} := (1_A \oplus 1_A)xx^*(1_A \oplus 1_A)$ has the desired
properties.

Now assume that $a \in \mathrm{M}_{\infty}(A)_{++}$.  Use Lemma \ref{specialrep}
to find representative $a^{'} \otimes 1_{\Z} \in A \otimes \Z$ of $\langle a \rangle$.
Put $b := a^{'} \otimes z_{1/r} \in \mathrm{M}_{\infty}(A \otimes \mathcal{Z})_+$,
so that $\dt(b) \leq 1$.  We now identify $A$ with $A \otimes \mathcal{Z}$ via $\phi$.
Our hypotheses ensure that $A$ has strict
comparison of positive elements (Corollary 4.6 of \cite{R1}),
whence $b \precsim 1_A$ by Proposition \ref{pureposcomparison}. We
apply~\cite[Proposition 2.4]{Rfunct} to $b + \epsilon\cdot 1_A
\precsim b \oplus \epsilon \precsim 1_A \oplus 1_A$, and obtain $x
\in \mathrm{M}_{\infty}(A)_+$ such that
\[
x^*(1_A \oplus 1_A)x = (b +\epsilon - \epsilon)_+ = b\,.
\]
It follows that
\[
b \sim \tilde{b} := (1_A \oplus 1_A)xx^*(1_A \oplus 1_A) \leq\Vert x\Vert^2 1_A \oplus 1_A\,.
\]
Now $(1/\Vert x\Vert^2) \tilde{b} \sim \tilde{b}$ --- Cuntz equivalence is robust
under multiplication by elements of $\mathbb{R}^{++}$ --- and so
\[
b \sim (1/ \Vert x \Vert^2) \tilde{b} \leq 1_A \oplus 1_A.
\]

It follows that
\[
(1/ \Vert x \Vert^2) (\tilde{b} \otimes z) \leq (1_A \oplus 1_A) \otimes
z,
\]
and that
\[
(1/\Vert x \Vert^2)(\tilde{b} \otimes z) \sim b \otimes z =(a^{'} \otimes z_{1/r}) \otimes
z
\]
(\cite[Lemma 4.1]{R1}). Put $\tilde{a} := (1/\Vert x \Vert^2)(\tilde{b} \otimes z)$.
The last equation shows that $\dt(\tilde{a}) = \dt(a)$, whence $a
\sim \tilde{a}$ by Proposition \ref{pureposcomparison}.
\end{proof}

\begin{props}\label{closure}
Let $A$ be a simple, unital, exact, and finite C$^*$-algebra absorbing the Jiang-Su algebra
$\mathcal{Z}$ tensorially.  Let there be given a sequence $(a_i)_{i=1}^{\infty} \subseteq
\mathrm{M}_{\infty}(A)_+$, and put
\[
 h_i(\tau) := \dt(a_i); \ \ g_i := \sum_{j=1}^{i} h_j.
\]
If
\[
\lim_{i \to \infty} g_i = g; \ \ \sum_{i=1}^{\infty} \Vert
h_i\Vert < \infty,
\]
then there exists $a \in \mathrm{M}_{\infty}(A)_{++}$ such that
$\dt(a) = g(\tau)$, for all $\tau \in \mathrm{T}(A)$.
\end{props}

\begin{proof}
We may assume that $a_i \in \mathrm{M}_{\infty}(A)_{++}$, since
$\dt(a_i) = \dt(a_i \otimes z_1)$, for all $\tau \in \mathrm{T}(A)$.
We may also assume that
$\sum_{i=1}^{\infty} \Vert h_i\Vert < 1$ by scaling the $a_i$
(using Corollary \ref{cone}).

Using the embedding of $C[0,1]$ into $\Z$ as in Theorem~\ref{WZ}
we may choose, for each $i \in \mathbb{N}$, a representative $y_i$
of $\langle z_{\Vert h_i\Vert} \rangle$ inside $\mathcal{Z}$ such
that $y_i y_j = y_j y_i = 0$ for all $i \neq j$. By
Lemma~\ref{specrep}, $a_i$ is equivalent to $\tilde{a}_i \leq
(1_A \oplus 1_A) \otimes y_i$. It follows that the $\tilde{a}_i$s are pairwise
orthogonal, and that $\dt(\tilde{a}_i) = h_i$. Put
\[
a := \sum_{i=1}^{\infty} \frac{1}{2^i} \tilde{a}_i \in
\mathrm{M}_2(A \otimes \mathcal{Z}).
\]
Then, $\dt(a) = g(\tau)$, as desired.
\end{proof}

Let $A$ be a C$^*$-algebra with finitely many pure tracial states.
In this situation we make the identifications
\[
\laff_b(\mathrm{T}(A))^{++}\equiv \aff(\mathrm{T}(A))^{++} \equiv
\{(\lambda_1,\ldots,\lambda_n)|\lambda_i \in \mathbb{R}^{++}\}\,,
\]
where $n$ is the number of pure tracial states on $A$.
Now suppose further that $A$ is simple, unital, exact, finite, and $\mathcal{Z}$-stable.
Since $\iota\colon W(A)_+\to\laff_b(\mathrm{T}(A))^{++}$ is an
order-embedding, we know (using~\cite[Theorem 2.6]{br}) that
$S((\mathbb{R}^{++})^n, 1)$ maps surjectively onto $S(W(A)_+,
\langle 1\otimes z_1\rangle)$, which by Lemma~\ref{k*struc} agrees
with
\[
S(\mathrm{K}_0^*(A),\mathrm{K}_0^*(A)^{++},[1_A]) \stackrel{\mathrm{def}}{=}\mathrm{DF}(A).
\]

\begin{rems}  {\rm
The definition of the term ``state'' is different for
partially ordered Abelian semigroups and partially ordered Abelian groups.
For semigroups a state must be order preserving, while for groups it is
required to be positive.  Both definitions require the state to be a linear map
into $\mathbb{R}$ taking the order unit to 1.  With this in mind it is easy to
check that the states on $(W(A),\langle 1_A \rangle)$ coincide with the states
on $(\mathrm{K}_0^*(A),\mathrm{K}_0^*(A)^{++},[1_A])$.  }
\end{rems}

Now, if $\tau$ is an extremal trace, then the corresponding lower
semicontinuous function $\dt$ is an extreme point in
$\mathrm{DF}(A)$. This follows from the fact that
$\mathrm{LDF}(A)$ is a face of $\mathrm{DF}(A)$ (\cite[Proposition
II.4.6]{bh}) and the fact that $\tau\mapsto\dt$ is an affine
bijection from $\mathrm{T}(A)$ onto $\mathrm{LDF}(A)$. In our case
of interest, where we have exactly $n$ extreme traces, we find
counting dimensions that
$S(\mathrm{K}_0^*(A),\mathrm{K}_0^*(A)^{++},[1_A])\cong
\mathbb{R}^n$. It follows from Corollary~\ref{k*archimedean}
and~\cite[Theorem 4.14]{G} that
$\mathrm{K}_0^*(A)\cong\mathbb{R}^n$ in this case.

Next, from the obvious containment
\[
S(\mathrm{K}_0^*(A),\mathrm{K}_0^*(A)^{++},[1_A])\subseteq
S(\mathrm{K}_0^*(A),\mathrm{K}_0^*(A)^{+}, [1_A])
\]
and the fact that $\mathrm{K}_0^*(A)\cong \mathbb{R}^n$, we see
that in fact we have equality.

We shall need the following result (see, e.g.~\cite[Theorem 7.9]{G}):

\begin{thms}\label{denseimage}
Let $(G,u)$ be an unperforated partially ordered Abelian group with order-unit, and let
\[
\psi\colon G \to \aff(S(G,u))
\]
be the natural map (given by evaluation). Then, the set
\[
\{\psi(x)/2^n\mid x \in G^+\,, \ n \in \mathbb{N}\}
\]
is dense in $\aff(S(G,u))^{+}$.
\end{thms}

Inspection of the proof reveals that the same result will hold under the assumption 
that $G$ is simple and weakly unperforated, which is what we shall use below.

\begin{thms}\label{ntraces}
Let $A$ be an exact, simple, and unital C$^*$-algebra absorbing
the Jiang-Su algebra $\mathcal{Z}$ tensorially.  Suppose that $A$
has $n$ pure tracial states.  Then, $\iota\colon W(A)_+ \to
\laff_b(\mathrm{T}(A))^{++}$ is surjective.
\end{thms}

\begin{proof}
From the comments preceding Theorem~\ref{denseimage}, it follows
that the state space of the group $\mathrm{K}_0^*(A)$ is
$\mathbb{R}^n$, no matter which ordering we consider on it (either
$\mathrm{K}_0^*(A)^+$ or $\mathrm{K}_0^*(A)^{++}$). Therefore,
\[
\aff(S(\mathrm{K}_0^*(A),\mathrm{K}_0^*(A)^+,[1_A]))=\aff(S(\mathrm{K}_0^*(A),\mathrm{K}_0^*(A)^{++},[1_A]))=\laff_b(\mathrm{T}(A))\,.
\]
We also know from Proposition~\ref{k*unperf} that
$(\mathrm{K}_0^*(A),\mathrm{K}_0^*(A)^+,[1_A])$ is a weakly unperforated
partially ordered simple abelian group. Our considerations above
together with Theorem~\ref{denseimage} imply that
\[
\{ \iota(a)/2^n\mid a \in \mathrm{M}_{\infty}(A)_{++}, \ n \in
\mathbb{N} \}
\]
is dense in $\laff_b(\mathrm{T}(A))$.  But $\iota(a)/2^n = \iota(a
\otimes z_{1/2^n})$ by Corollary \ref{cone}, so
\[
\{ \iota(a)/2^n\mid  a \in \mathrm{M}_{\infty}(A)_{++}, \ n \in
\mathbb{N} \} = \{ \iota(a)\mid  a \in \mathrm{M}_{\infty}(A)_{++}
\}.
\]
In other words, the image of $\iota$ in
$\laff_b(\mathrm{T}(A))^{++}$ is dense.

Let $f \in \laff_b(\mathrm{T}(A))^{++}$ be given. A moment's
reflection shows that one may choose a sequence
$(h_i)_{i=1}^{\infty} \subseteq \laff_b(\mathrm{T}(A))^{++}$ with
the following properties:
\begin{enumerate}
\item[(i)] $\lim_{i \to \infty} f_i = f$, where $f_i =
\sum_{j=1}^i h_j$; \item[(ii)] $\sum_{i=1}^{\infty} \Vert h_i\Vert
< \infty$; \item[(iii)] $h_i(\tau) = \dt(a_i)$ for some $a_i \in
\mathrm{M}_{\infty}(A)_{++}$.
\end{enumerate}

We may apply Proposition \ref{closure} to find $a \in
\mathrm{M}_{\infty}(A)_{++}$ such that $\dt(a) = f(\tau)$, for all
$\tau \in \mathrm{T}(A)$, whence $\iota$ is surjective, as
desired.
\end{proof}

\section{Real rank zero}

In this section we show that our map $\iota$ is surjective
whenever $A$ is a $\Z$-stable, simple, exact C$^*$-algebra with
real rank zero and stable rank one. In fact, we can prove a more
general result, namely that for such an $A$ (not necessarily
simple) $\mathrm{K}_0^*(A)$ is order-isomorphic to the group of
differences of lower semicontinous, affine, real-valued and
bounded functions defined on $\mathrm{T}(A)$, equipped with the
pointwise ordering. Some of our arguments, namely the first part
of Theorem~\ref{rr0} below, can be traced back to the ones
in~\cite{bh}, and we include them for the convenience of the
reader.

It should be no surprise, however, that the (WEC) implies the (EC)
for this class. This can be justified by recalling that the Cuntz
semigroup $W(A)$ is completely determined by $V(A)$ whenever $A$
is $\sigma$-unital, has real rank zero and stable rank one. More
concretely, one can obtain for such an $A$ an order-isomorphism
between $W(A)$ and the monoid of the so-called countably generated
intervals in $V(A)$ that are bounded by the generating interval
$D(A)$ (see~\cite{P1} for a full account).

Given a positively ordered abelian semigroup with order-unit
$(M,\leq, u)$, consider the natural representation map
$\phi_u\colon M\to \aff(S(M,u))^+$. It is said that $M$ satisfies
\emph{condition (D)} provided that $\phi_u(M)$ is dense. A unital
C$^*$-algebra $A$ satisfies condition (D) provided that the
positive cone $\mathrm{K}_0(A)^+$ of its Grothendieck group
satisfies condition (D). It was shown in~\cite{pardo} that any
unital C$^*$-algebra $A$ with real rank zero satisfies condition
(D) if and only if $A$ has no finite dimensional representations.

\begin{lms}
\label{almost} Let $A$ be a $\Z$-stable unital C$^*$-algebra with
stable rank one. Then $s(x)>0$ for all states $s$ on
$S(\mathrm{K}_0(A),[1_A])$ if and only if $x$ is an order-unit in
$\mathrm{K}_0(A)$.
\end{lms}

\begin{proof}
Since $A$ has stable rank one, we have $\mathrm{K}_0(A)^+=V(A)$.
We also know from~\cite[Corollary 4.8]{R1} that $V(A)$ is almost
unperforated. Assume that $s(x)>0$ for all states $s$. It then
follows from~\cite[Theorem 4.12]{G} that $mx$ is an order-unit for
some natural number $m$. Write $x=a-b$ where $a$, $b\in V(A)$. We
know that there is $l$ in $\mathbb{N}$ such that $b\leq lm(a-b)$,
and hence $(lm+1)b\leq lma$. Therefore $b\leq a$, and so $x>0$.
Thus $x$ is an order-unit.
\end{proof}

If $f$, $g$ are real-valued functions defined on a set $X$, write
$f\gg g$ (or $f\ll g$) to mean that $f(x)>g(x)$ (or $f(x)<g(x)$) for every $x$ in $X$.

\begin{lms}
\label{orth} Let $A$ be a $\Z$-stable unital C$^*$-algebra with
real rank zero and stable rank one. Then $A$ contains a sequence
of orthogonal projections $(p_n)$ such that $s([p_n])>0$ for all
states $s\in S(V(A),[1_A])$. (Equivalently, $\tau(p_n)>0$ for all
quasitraces on $A$.)
\end{lms}

\begin{proof} (Outline.) Note first that $A\cong A\otimes\Z$ satisfies condition
(D), because $\Z$ is simple and infinite dimensional. Denote by
$u=[1_A]\in V(A)$ and by
\[
\phi_u\colon V(A)\to \aff(S(V(A),u))=\aff(S(\mathrm{K}_0(A),u))
\]
the natural representation map, given by evaluation.

Using condition (D) we may then find a projection $p_1$ such that
$0\ll\phi_u([p_1])\ll 1$. Thus, by compactness of the state space
of $V(A)$ and condition (D) again, there is a projection $p_2'$
satisfying $0\ll\phi_u([p_2'])\ll\phi_u([1-p_1])$.
Lemma~\ref{almost} implies that $p_2'\sim p_2\leq 1-p_1$ for some
projection $p_2$. Continuing in this way we find our sequence of
projections $(p_n)$.

The equivalent statement follows readily from the fact that the
map $\mathrm{QT}(A)\to S(V(A),[1_A])$, given by evaluation, is an
affine homeomorphism (see~\cite[Theorem III.1.3]{bh}).
\end{proof}

We remark that Lemma~\ref{orth} also holds trading $\Z$-stability and stable rank one by weak divisibility. This latter property was introduced in~\cite{peror}: a C$^*$-algebra $A$ is \emph{weakly
divisible} if for any element $x$ in $V(A)$, we may find a natural number $n$ and elements $y$ and $z$ in
$V(A)$ such that $x=ny+(n+1)z$. Weak divisibility is always
guaranteed for simple (non-type I) C$^*$-algebras of real rank
zero, and holds quite widely in the non-simple case
(see~\cite[Theorem 5.8]{peror}). Basically, what we need to use to establish~\ref{orth} in this setting is that for a non-zero $x$ in $V(A)$, there is $n$ and a non-zero $y$ in $V(A)$ such that $ny\leq x\leq (n+1)y$.

\begin{thms}
\label{rr0} {\rm (cf. \cite[Theorem III.3.2 and Corollary III.3.3]{bh})} Let $A$ be a $\Z$-stable,
exact, separable and unital C$^*$-algebra with real rank zero and stable rank one. Then
$\mathrm{K}_0^*(A)$ is order-isomorphic to $G(\laff_b(\mathrm{T}(A))$, equipped with the pointwise ordering.
\end{thms}

\begin{proof}
Define $\iota\colon \mathrm{K}_0^*(A)\to
G(\laff_b(\mathrm{T}(A)))$ by $\iota([a])(\tau)=\dt(a)$. Note
first that, for a positive element $a$, if $(p_n)$ is an
(increasing) approximate unit consisting of projections for the
hereditary algebra generated by $a$, we have that
$\iota([a])(\tau)=\sup_n\tau(p_n)$.

In order to get an order-isomorphism onto the image, we have to show that $[a]\leq [b]$ in
$\mathrm{K}_0^*(A)$ whenever $\iota([a])\leq\iota([b])$. Let $(p_n)$ be the sequence of
orthogonal projections constructed in Lemma~\ref{orth}, and let $c=\sum_{n=1}^{\infty}\frac{1}{2^n}r_n\in A_+$,
where $r_n=\sum_{i=1}^n p_i$. Let $(e_n)$ and $(f_n)$ be approximate units consisting of projections
for the hereditary algebras generated by $a$ and $b$ respectively. We then have that $(e_n\oplus r_n)$
(respectively, $(f_n\oplus r_n)$) is an (increasing) approximate unit consisting of projections for
$a\oplus c$ (respectively, for $b\oplus c$). Note that $\iota([a\oplus c])\leq\iota([b\oplus c])$.
By construction of the sequence $(r_n)$ and Lemma~\ref{orth}, the sequence $\tau(e_n\oplus r_n)$ is
strictly increasing. Using compactness of the state space of $V(A)$, we find that for all $n$, there
is $m$ such that $\tau(e_n\oplus r_n)<\tau(f_m\oplus r_m)$ for all $\tau$. It follows again from
Lemma~\ref{orth} that for all $n$, there is $m$ such that  $e_n\oplus r_n\precsim f_m\oplus r_m$.
But this implies that $a\oplus c\precsim b\oplus c$ (see~\cite[Proposition 2.3]{P1} and also~\cite[Corollary III.3.8]{bh}).

We now prove that $\iota$ is surjective. Let $f\in\laff_b(\mathrm{T}(A))$, which is bounded below
by some constant $k$. Writing $h=f-k+1$, we may assume that actually $f\in\laff_b(\mathrm{T}(A))^{++}$. Since $A$ is separable, we have
that $\mathrm{T}(A)$ is metrizable, hence we may write $f$ as a
pointwise supremum of an increasing sequence $(f_n)$ of functions
in $\aff (\mathrm{T}(A))^{++}$. Choose $n_0$ such that
$f_n-\frac{1}{2^n}\gg 0$ whenever $n\geq n_0$. Write $u=[1_A]\in V(A)$ and denote as before $\phi_u$
the natural representation map.

Using condition (D) we may find projections $p_n$ in $M_{\infty}(A)$ such that
$f_n-\frac{1}{2^n}\ll\phi_u([p_n])\ll f_n-\frac{1}{2^{n+1}}$ for
all $n\geq n_0$, where $u=[1_A]\in V(A)$. Since
$\phi_u([p_n])\ll\phi_u([p_{n+1}])$ we get from Lemma~\ref{almost}
that $[p_n]\leq [p_{n+1}]$ in $V(A)$. Since $f$ is also bounded, a
second use of Lemma~\ref{almost} shows that $p_n$ all belong to
$M_t(A)$ for some $t$. Using that $A$ has stable rank one (whence
projections cancel from direct sums) we may arrange that the
sequence $(p_n)$ is indeed increasing in the order of $A$.

It is clear that $f$, being the pointwise supremum of the $f_n$'s,
will satisfy that $f=\sup \phi_u([p_n])$. We know from~\cite[Theorem III.1.3]{bh} that the
natural mapping $\mathrm{T}(A)\to S(\mathrm{K}_0(A),[1_A])$ is a homeomorphism.

If we then let $x=\sum\limits_{n=1}^{\infty}\frac{1}{2^n}p_n$, we find that $x\otimes z_1$ is
a purely positive element in $M_t(A)$ such that $\dt(x\otimes z_1)=\dt(x)=\sup_n\dt(p_n)=\sup
\tau(p_n)=\phi_u([p_n])(\tau)=f(\tau)$ for every $\tau\in\mathrm{T}(A)$.
\end{proof}

The argument of surjectivity in the proof of Theorem~\ref{rr0}, allows us to state the following:

\begin{cors}
\label{rr0simple}
Let $A$ be an exact, simple, and unital C$^*$-algebra absorbing
the Jiang-Su algebra $\mathcal{Z}$ tensorially. Suppose that $A$ has real rank zero and stable rank one.
Then, $\iota\colon W(A)_+ \to \laff_b(\mathrm{T}(A))^{++}$ is surjective.
\end{cors}

\section{Goodearl algebras}

In this section we prove that $\iota$ is surjective for algebras we term
{\it degenerate Goodearl algebras}, and outline a proof of the same fact
for the simple Goodearl algebras studied in \cite{gpubmat}.
We do so to support the conjecture that $\iota$ is always surjective for
unital and stably finite C$^*$-algebras without nonzero finite-dimensional
representations. 
In other words, hypotheses (i) and (ii) of Theorem \ref{equiv} should be
removeable.  (Note that for a non-simple algebra, the image of $\iota$
will not always consist of strictly positive functions.)  
 Our reasons for providing a sketch in lieu of a full
proof in the simple case are twofold:  first,
the main ideas and technical details for a full proof are contained 
already in the argument for the degenerate case;  second, simple 
Goodearl algebras are known to satisfy the Elliott conjecture, 
and so one gains little new insight into their structure by computing 
their Cuntz semigroups.

Let $X$ and $Y$ be compact Hausdorff spaces. A $*$-homomorphism
\[
\phi\colon \co(X) \to \ma_n(\co(Y))
\]
is called \emph{diagonal} if
\[
\phi(f)(y) = \mathrm{diag}\left(f(\gamma_1(y)),\ldots,f(\gamma_n(y))\right)
\]
for continuous maps $\gamma_i\colon Y \to X$, $1 \leq i \leq n$.  The $\gamma_i$
are called \emph{eigenvalue maps}.

Let $X$ be a nonempty, separable, and compact Hausdorff space. Let $A =
\lim_{i \to \infty}(A_i,\phi_i)$ be a unital inductive limit
C$^*$-algebra where, for each $i \in \mathbb{N}$, $A_i \cong
\ma_{n_i}(\co(X))$ for some $n_i \in \mathbb{N}$ with $n_i |
n_{i+1}$, $\phi_i$ is diagonal, and the eigenvalue maps of
$\phi_i$ are either the identity map on $X$, or have range equal
to a point. Such an algebra will be called a \emph{Goodearl
algebra}.  This definition generalises slightly the one provided
by Goodearl in \cite{gpubmat}.

If each $\phi_i$ in the inductive sequence for $A$ has every
eigenvalue map equal to the identity map on $X$, then we will say
that $A$ is \emph{degenerate}. In this case one obtains a
(in general non-simple) algebra isomorphic to the tensor product $\co(X)
\otimes \mathfrak{U}$, where $\mathfrak{U}$ is the UHF algebra
whose $\mathrm{K}_0$-group is the subgroup of the rationals whose
denominators, when in lowest terms, divide some $n_i$.  This
subgroup is dense in $\mathbb{R}$ whenever $n_i \to \infty$ as $i
\to \infty$. In this case, $\mathrm{T}(A)$ may be identified with the Bauer
simplex $M_1^+(X)$ of positive probability measures on $X$, hence its extreme boundary $\partial_e \mathrm{T}(A)$
is homeomorphic to $X$.  (Recall that a Bauer simplex is a Choquet simplex
with closed extreme boundary --- see \cite{Al} for details.)

If $X$ is a compact Hausdorff space, denote by $\mathrm{L}(X)$ the
semigroup of lower semicontinuous real-valued functions defined on
$X$, by $\mathrm{L}(X)^{++}$ the subsemigroup of $\mathrm{L}(X)$
consisting of strictly positive elements, and by $\mathrm{L}_b(X)$
the subsemigroup of bounded functions. Let $A$ be a unital
C$^*$-algebra such that $\mathrm{T}(A)$ is a non-empty Bauer
simplex.  Then, there is a semigroup isomorphism between
$\laff_b(\mathrm{T}(A))$ and $\mathrm{L}_b(\partial_e
\mathrm{T}(A))$ -- the behaviour of $f \in \laff(\mathrm{T}(A))$
is determined by the behaviour of its restriction to $\partial_e
\mathrm{T}(A)$ (cf.~\cite[Lemma 7.2]{GK}). It follows that proving the surjectivity of $\iota$
for such an algebra only requires proving that every $f \in
\mathrm{L}_b(\partial_e\mathrm{T}(A))^{++}$ can be realised as the
image of some $a \in \mathrm{M}_{\infty}(A)_{++}$ under the map
\[
\iota_e\colon  W(A)_+ \to \mathrm{L}_b(\partial_e\mathrm{T}(A))^{++}
\]
given by
\[
\iota_e(\langle a\rangle ) = \dt(a), \ \text{for all } \tau \in
\partial_e \mathrm{T}(A).
\]
Clearly, it will suffice to prove the above for functions $f$ such
that $\Vert f\Vert \leq 1$.

\begin{thms}
\label{nonsimple} Let $A$ be a degenerate Goodearl algebra. Then
$\iota$ is surjective.
\end{thms}
\begin{proof}

We identify $\mathrm{T}(A)$ with the Bauer simplex $M_1^+(X)$,
whence $\partial_e(\mathrm{T}(A))$ is homeomorphic to $X$. Let us write $\tau_x$
for the trace that corresponds to a point $x$ in $X$. This, in turn, corresponds to the point
mass measure $\delta_x$ at $x$.

Let $f \in \mathrm{L}_b(X)^{++}$ be given, and assume that $\Vert f\Vert \leq 1$.  We
prove that $f$ is the image of an element $a \in A_+$ under the map
$\iota_e$ defined above.

Define, for each $i \in \mathbb{N}$, a function $f_i$
as follows:  put
\[
F_{i,k} := \left\{x \in X\mid f(x) \leq \frac{k}{n_i} \right\}\,, \ 1
\leq k \leq n_i\,,
\]
\[f_i(x) = 0\,, \text{for all } x \in F_{1,k}\,,
\]
and
\[
f_i(x) := \frac{k-1}{n_i} \ \text{whenever }  x \in F_{i,k}
\setminus F_{i,k-1}\,.
\]
Let us check that $f_i$ converges pointwise to $f$, and that $f_j \geq f_i$
whenever $j \geq i$. Let $x\in F_{i,k}\setminus F_{i,k-1}$ for $1\leq k\leq n_i$, and
take $j\geq i$. Then $f_i(x)=\frac{k-1}{n_i}$. Write $n_j=n_in_i'$, and note that
$f(x)\leq\frac{k}{n_i}=\frac{kn_i'}{n_j}$. Thus $x\in F_{j,kn_i'}$. Let $l\geq 0$
be such that $x\in F_{j,kn_i'-l}\setminus F_{j,kn_i'-l-1}$. Since
\[
\frac{k-1}{n_i}<f(x)\leq\frac{kn_i'-l}{n_j}\,,
\]
it is easy to check now that $f_j(x)=\frac{kn_i'-l-1}{n_j}\geq\frac{k-1}{n_i}=f_i(x)$.

Note that for $x\in F_{i,k}\setminus F_{i,k-1}$ we have $f(x)-f_i(x)\leq \frac{1}{n_i}$,
whence clearly $f_i\rightarrow f$.

We will construct an increasing sequence $a_1 \leq a_2 \leq \ldots$ of positive elements
in $A$ converging to a
positive element $a$, such that $d_{\tau_x}(a_i) = f_i(x)$, for all $x\in X$. It will
follow that $d_{\tau_x}(a) =f(x)$, for all $x\in X$.

For each $i \in \mathbb{N}$, choose $n_i$ continuous functions
$f_{i,k}\colon X \to [0,1/2^i]$ as follows: $f_{i,1} \equiv 0$,
and $f_{i,k}$ is supported on the open set $F_{i,k-1}^{c}$, for $2
\leq k \leq n_i$. Put
\[
\tilde{a}_i := \mathrm{diag}(f_{i,1},\ldots,f_{i,n_i}) \in A_i\,.
\]

Define $a_1:=\tilde{a}_1\in A_1$. Suppose that we have constructed $a_1,\ldots,a_i$
such that $a_j\in A_j$ and also $a_1\leq a_2\leq\ldots\leq a_i$ when viewed in $A_i$
(through the natural maps). We now construct $a_{i+1}$.

Consider the image of $a_i$ in $A_{i+1}$ under $\phi_i$.  It is a diagonal
element, and its diagonal entries consist of $n_{i+1}/n_i$ copies
of $f_{i,k}$ for each $1 \leq k \leq n_i$.
Now, for any such $k$, notice that the open set $F_{i,k}^{c}$ is contained in $F_{i+1,l}^{c}$
for every $(k-1)(n_{i+1}/n_i)+1 \leq l \leq k(n_{i+1}/n_i)$.
Assume, by permuting the diagonal entries of $\tilde{a}_{i+1}$ if
necessary, that the entries of $\phi(a_i)$ equal to $f_{i,k}$
correspond to the entries of $f_{i+1,l}$ of $\tilde{a}_{i+1}$ for
which $(k-1)(n_{i+1}/n_i)+1 \leq l \leq k(n_{i+1}/n_i)$. Now
define $a_{i+1}$ to be the diagonal element whose entries are the
pointwise maximum of the entries of $\phi_i(a_i)$ and
$\tilde{a}_{i+1}$.

Since $F_{i,k}^{c}\subseteq F_{i+1,l}^{c}$, we have that
$\mathrm{Coz}(\max\{f_{i,k}, f_{i+1,l}\})=\mathrm{Coz}(f_{i+1,l})=F_{i+1,l}^{c}$
($\mathrm{Coz}(f)$ denotes the cozero set of a function $f$).
For any $x\in X$, we have
\[
d_{\tau_x}(a_{i+1})=d_{\tau_x}(\tilde{a}_{i+1})=\frac{1}{n_{i+1}}\sum_{j=1}^{n_{i+1}}\delta_x(F_{i+1,j}^{c})=\frac{k}{n_{i+1}}\,,
\]
where $k$ is such that $x\in F_{i+1,k}^{c}\setminus F_{i+1,k+1}^{c}$.  Hence $d_{\tau_x}(a_{i+1})=f_{i+1}(x)$.
Observe that $\phi_i(a_i)\leq a_{i+1}$ and $\Vert a_i- a_{i-1}\Vert<1/2^i$ by construction.

Continue in this way and identify the $a_i$'s with their images in $A$. Then the sequence $(a_i)_{i=1}^{\infty} \subseteq A$ has the following
properties:
\begin{enumerate}
\item[(i)] $a_i \leq a_{i+1}$ for all $i$; \item[(ii)] $\Vert a_i-a_{i-1}\Vert < 1/2^i$;
\item[(iii)] $d_{\tau_x}(a_i) = f_i(x)$, for all
$x\in X$.
\end{enumerate}
It follows that $a:=\lim_{j \to \infty} a_j$ has the desired
property:
\[
d_{\tau_x}(a) = f(x), \ \text{for all } x\in X.
\]
\end{proof}

Simple Goodearl algebras are either of real rank zero or real rank one (\cite{gpubmat}),
and are known to be approximately divisible (see \cite{EGL2}).  It follows from Theorem 2.3
of \cite{TW2} that they are $\mathcal{Z}$-stable, and so, in the real rank zero case, the surjectivity
of $\iota$ is given by Corollary~\ref{rr0simple}.  In the real rank one case, it
is known that the connecting $*$-homomorphisms $\phi_i$ in the inductive sequence for
the given algebra must contain a vanishingly small proportion of eigenvalue maps with
range equal to a point --- the connecting maps are very nearly those of a degenerate Goodearl
algebra (\cite{gpubmat}).  Combining this fact with the construction of Theorem~\ref{nonsimple}
yields the surjectivity of $\iota$ for simple Goodearl algebras of real rank one.
The details are left to the reader.

\section{Concluding remarks}

Although $\mathcal{Z}$-stability is a useful tool in the proofs of Theorem \ref{ntraces}
and Corollary \ref{rr0simple}, it is by no means a necessary condition for the surjectivity of
$\iota$.  A calculation akin to the proof of Theorem \ref{nonsimple} shows that $\iota$ is surjective for
the non-$\mathcal{Z}$-stable AH algebra constructed in Theorem 1.1 of \cite{T2}.  Also:

\begin{props}\label{onetracesurj}
Let $A$ be a unital C$^*$-algebra with unique tracial state $\tau$.  Suppose that there exists $a \in A^+$
such that $\mathrm{Sp}(a) = [0,1]$, and that $\tau$ induces an atom-free measure on $\mathrm{Sp}(a)$.
Then, $\iota$ is surjective.
\end{props}

\begin{proof}
We need only produce, for every $\lambda \in (0,1]$, positive elements $g_{\lambda} \in A$
such that $\dt(g_{\lambda}) = \lambda$.  This is straightforward:  let $O_{\lambda}$ be
an open set of measure $\lambda$ with respect to $\tau$ (such a set exists since said
measure is an atom-free probability measure on $[0,1]$), and let $g_{\lambda}$ be a
positive function supported on $O_{\lambda}$.
\end{proof}

The results of Sections 6, 7, and 8 suggest a closing question:
\begin{qus}\label{surj}
Is $\iota$ surjective for any unital and stably finite C$^*$-algebra $A$ having no
nonzero finite-dimensional representations?
\end{qus}

\noindent
An affirmative answer will extend the equivalence of (EC) and (WEC) to all
simple, separable, unital, nuclear, finite, and $\mathcal{Z}$-stable C$^*$-algebras.


\markboth{}{}

\begin{thebibliography}{99}

\bibitem{Al} Alfsen, E.: Compact convex sets and boundary integrals.  Ergebnisse der Mathematik und ihrer
Grenzgebiete, Band 57, Springer-Verlag, 1971.

\bibitem{blsur} Blackadar, B.: \emph{Comparison theory for simple C$^*$-algebras}, in Operator Algebras and Applications, eds.
D. E. Evans and M. Takesaki, LMS Lecture Notes Series, \textbf{135}, Cambridge Univ. Press, 1988, pp. 21--54.

\bibitem{bh} Blackadar, B.\ and Handelman, D.: \emph{Dimension functions and traces on C$^*$-algebras}, 
J. Funct. Anal. \textbf{45} (1982), pp. 297--340.

\bibitem{br} Blackadar, B.\ and R\o rdam, M.: \emph{Extending states on
preordered semigroups and the existence of quasitraces on
C$^*$-algebras}, J. Algebra \textbf{152} (1992), pp. 240--247.

\bibitem{Cu} J. Cuntz: \emph{Dimension functions on simple C$^*$-algebras}, Math. Ann. \textbf{233} (1978), pp. 145--153.

\bibitem{DG} Dadarlat, M.\ and Gong., G.: \emph{A classification result for
approximately homogeneous C$^*$-algebras of real rank zero}, Geom.\ Funct.\
Anal.\ {\bf 7} (1997), pp. 646-711.

\bibitem{Ei} Eilers, S.: \emph{A complete invariant for $AD$ algebras of real rank zero with bounded torsion in $\mathrm{K}_1$},
J. Funct. Anal. {\bf 139} (1996), pp. 325-348.

\bibitem{El1} Elliott, G.\ A.: \emph{On the classification of inductive
limits of sequences of semi-simple finite-dimensional algebras}, J.\ Algebra
{\bf 38} (1976), pp. 29-44.

\bibitem{El2} Elliott, G.\ A.: \emph{The classification problem for amenable
C$^*$-algebras}, Proc. ICM '94, Zurich, Switzerland, Birkhauser Verlag, Basel,
Switzerland, pp. 922-932.

\bibitem{El3} Elliott, G.\ A.: \emph{On the classification of C$^*$-algebras
of real rank zero}, J.\ Reine Angew.\ Math.\ {\bf 443} (1993), pp. 179-219.

\bibitem{EG} Elliott, G. A.\ and Gong, G.: \emph{On the classification of C$^*$-algebras of real rank zero. II.}, 
Ann. of Math. (2) {\bf 144} (1996), pp. 497-610.

\bibitem{EGL} Elliott, G.\ A., Gong, G.\ and Li, L.: \emph{On the classification
of simple inductive limit C$^*$-algebras, II:  The isomorphism theorem}, preprint.

\bibitem{EGL2} Elliott, G.\ A., Gong, G.\ and Li, L.: \emph{Approximate divisibility of simple inductive limit 
C$^*$-algebras}, Operator algebras and operator
theory (Shanghai, 1997), 87-97, Comtemp. Math. 228.

\bibitem{gjsu} Gong, G., Jiang, X.\ and Su, H.: \emph{Obstructions to $\Z$-stability for unital simple
C$^*$-algebras}, Canad. Math. Bull. \textbf{43} (2000), pp.
418--426.

\bibitem{G} Goodearl, K.\ R.: \emph{Partially Ordered Abelian Groups with Interpolation}, Math. 
Surveys and Monographs \textbf{20}, Amer. Math. Soc., Providence, 1986.

\bibitem{gpubmat} Goodearl, K.\ R.: \emph{Notes on a class of simple C$^*$-algebras with real rank
zero}, Publ. Mat. \textbf{36} (1992), pp. 637--654.

\bibitem{GK} Goodearl, K.\ R.: \emph{$\mathrm{K}_0$ of multiplier algebras of C$^*$-algebras with real rank
zero}, K-Theory \textbf{10} (1996), pp. 419--489.

\bibitem{Ha} Haagerup, U.: \emph{Quasitraces on exact C$^*$-algebras are traces}, preprint, 1991.

\bibitem{han} Handelman, D.: \emph{Homomorphisms of C$^*$-algebras to finite $AW^*$-algebras}, Michigan Math. J. \textbf{28} (1981), pp. 229-240.

\bibitem{JS1} Jiang, X.\ and Su, H.: \emph{On a simple unital projectionless $C^{*}$-algebra}, Amer.\ J.\ Math.\ {\bf 121} (1999), pp. 359-413.

\bibitem{K} Kirchberg, E.: \emph{The classification of Purely Infinite C$^*$-algebras using Kasparov's Theory}, in preparation.

\bibitem{KR} Kirchberg, E.\ and R\o rdam, M.: \emph{Non-simple purely infinite C$^*$-algebras}, Amer. J. Math. \textbf{122} (2000), pp.
637--666.

\bibitem{Li1} Lin, H.:  \emph{Classification of simple tracially AF
$C^{*}$-algebras},  Canad.\ J.\ Math.\  {\bf 53}  (2001),  pp. 161-194.

\bibitem{LZ} Lin, H.\ and Zhang, S.: \emph{On infinite simple
C$^*$-algebras}, J. Funct. Anal. \textbf{100} (1991), pp.
221-231.

\bibitem{M} Mygind, J.: \emph{Classification of certain simple C$^*$-algebras with torsion in $\mathrm{K}_1$},
Canad.\ J.\ Math.\ {\bf 53} (2001), pp. 1223-1308.
\bibitem{pardo} E. Pardo, \emph{Metric completions of ordered groups and $\mathrm{K}_0$ of exchange rings}, 
Trans. Amer. Math. Soc. \textbf{350} (1998), pp. 913-933.

\bibitem{P1} Perera, F.: \emph{The structure of positive elements for C$^*$-algebras with real rank zero}, 
International J. Math. \textbf{8} (1997), pp. 383-405.

\bibitem{percan} Perera, F.: \emph{Ideal structure of multiplier algebras of simple C$^*$-algebras with real rank
zero}, Canad. J. Math. \textbf{53} (2001), pp. 592-630.

\bibitem{peror} Perera, F.\ and R\o rdam, M.: \emph{$AF$-embeddings into
C$^*$-algebras of real rank zero}, J. Funct. Anal. \textbf{217}
(2004), pp. 142-170.
\bibitem{P} Phillips, N.\ C.: \emph{A classification theorem for nuclear
purely infinte simple C$^*$-algebras}, Doc.\ Math.\ {\bf 5} (2000), pp. 49-114.

\bibitem{Rfunct} R\o rdam, M.: \emph{On the structure of simple C$^*$-algebras tensored with a UHF-algebra. II}, 
J. Funct. Anal. \textbf{107} (1992), pp. 255-269.

\bibitem{R2} R\o rdam, M.: \emph{A simple C$^*$-algebra with a finite and an infinite projection}, 
Acta Math. \textbf{191} (2003), pp. 109-142.

\bibitem{R1} R\o rdam, M.: \emph{The stable and the real rank of \Z-absorbing C$^*$-algebras}, 
International J. Math. \textbf{15} (2004), pp. 1065-1084.

\bibitem{R3} R\o rdam, M.: \emph{Classification of Nuclear C$^*$-Algebras}, Encyclopaedia of Mathematical Sciences {\bf 126}, Springer-Verlag, 2002.

\bibitem{Th} Thomsen, K.: \emph{Limits of certain subhomogeneous C$^*$-algebras}, Mem. Soc. Math. Fr. \textbf{71} (1997), vi+125 pp. (1998).

\bibitem{T1} Toms, A.\ S.: \emph{On the independence of $\mathrm{K}$-theory and stable rank for 
simple C$^*$-algebras}, J. Reine und Angew. Math. {\bf 578} (2005), pp. 185-199.

\bibitem{T2} Toms, A.\ S.: \emph{On the classification problem for nuclear C$^*$-algebras}, arXiv preprint math.OA/0509103 (2005).

\bibitem{TW1} Toms, A.\ S.\ and Winter, W.: {\it Strongly self-absorbing C$^*$-algebras}, to appear in Trans. Amer. Math. Soc.

\bibitem{TW2} Toms, A.\ S.\ and Winter, W.: \emph{$\mathcal{Z}$-stable ASH algebras}, to appear in Canad. J. Math.

\bibitem{V1} Villadsen, J.: \emph{Simple C$^*$-algebras with perforation}, J.\ Funct.\ Anal.\ {\bf 154} (1998), pp. 110-116.

\bibitem{V2} Villadsen, J.: \emph{On the stable rank of simple C$^*$-algebras}, J.\ Amer.\ Math.\ Soc.\ {\bf 12} (1999), pp. 1091-1102.

\end{thebibliography}
\end{document}